\documentclass[12pt,a4paper]{amsart}
\usepackage[latin1]{inputenc}
\usepackage[T1]{fontenc}
\usepackage[french]{babel}
\usepackage{amssymb,amsfonts,latexsym}
\usepackage[arrow,matrix,curve]{xy}
\usepackage[mathscr]{eucal}
\usepackage{a4wide}
\usepackage{enumerate}

\font\trmbi=cmtcsc10 scaled 1200
\font\srmbi=cmtcsc10 scaled 840         
\font\ssrmbi=cmtcsc10 scaled 600
\font\teuf=eufm10 scaled 1200
\font\seuf=eufm7 scaled 1200            
\font\sseuf=eufm5 scaled 1200
\font\tenscr=rsfs10  scaled 1200
\font\sevenscr=rsfs7 scaled 1200        
\font\fivescr=rsfs5  scaled 1200

\skewchar\tenscr='177 \skewchar\sevenscr='177 \skewchar\fivescr='177

\newfam\majfam
        \textfont\majfam=\trmbi
        \scriptfont\majfam=\srmbi
        \scriptscriptfont\majfam=\ssrmbi

\newfam\euffam
        \textfont\euffam=\teuf
        \scriptfont\euffam=\seuf
        \scriptscriptfont\euffam=\sseuf

\newfam\scrfam 
        \textfont\scrfam=\tenscr 
        \scriptfont\scrfam=\sevenscr
        \scriptscriptfont\scrfam=\fivescr

\def \goth{\fam\euffam}

\def \scr{\fam\scrfam}

\def \SSS{\scriptscriptstyle}

\def \iso{\hbox{$\,\vbox{\hbox to
18pt{\hfill$\sim$\hfill}\kern-11pt\hbox{$\longrightarrow$}}\,$}}
\def \liso{\hbox{$\,\vbox{\hbox to
18pt{\hfill$\sim$\hfill}\kern-11pt\hbox{$\longleftarrow$}}\,$}}
\def \doublefleche{\vbox
to5.5025pt{\hbox{$\longrightarrow$}\kern-10.75pt\hbox{$\longrightarrow$}}}
\def \subsetneq{\;\lower1.66667pt\vbox
to9.60166pt{\vskip0.3pt\hbox{$\subseteq$}\kern-12.83pt\hbox{$\;\SSS/\;$}}\;}

\def \subsetneq{\;\lower1.66667pt\vbox
to9.60166pt{\vskip0.3pt\hbox{$\subseteq$}\kern-12.83pt\hbox{$\;\SSS/\;$}}\;}
\def \iso{\buildrel\sim\over\longrightarrow}

\def \HM{{\mathbb{H}}}

\def \NM{{\mathbb{N}}}

\def \PM{{\mathbb{P}}}
\def \QM{{\mathbb{Q}}}
\def \RM{{\mathbb{R}}}

\def \ZM{{\mathbb{Z}}}

\def \ACC{{\scr A}}

\def \ECC{{\scr E}}
\def \FCC{{\scr F}}

\def \ICC{{\scr I}}
\def \JCC{{\scr J}}
\def \KCC{{\scr K}}
\def \LCC{{\scr L}}

\def \OCC{{\scr O}}

\def \UCC{{\scr U}}
\def \VCC{{\scr V}}

\def \YCC{{\scr Y}}

\def \PG{{\goth P}}

\def \UG{{\goth U}}

\def \TS{{\mathsf T}}

\def \YS{{\mathsf Y}}

\def \leq{\leqslant}
\def \geq{\geqslant}

\newcommand{\rqe}{\noindent {\sc Remarque : }\rm}
\newcommand{\notation}{\noindent {\sc Notations : }\rm}

\newcommand {\dem}{\noindent {\sc D{\'e}monstration : } \rm }
\newcommand {\findem}{\hfill$\Box$\par\medskip}
\newcount\exno
\newcommand{\exo}[1]{\advance\exno by1 \medskip
                \hangindent=\parindent\hangafter=1
                \noindent{\trmbi Exercice}
                \hbox{ \arabic{chapter}.\the\exno.}\hfill \\ 
               \rm #1\par\medskip}

\renewcommand{\thesubsection}{\arabic{section}.\arabic{subsection}}

\newtheorem{prop}{Proposition}[subsection]
\newtheorem{lemme}[prop]{Lemme}
\newtheorem{Defi}[prop]{D{\'e}finition}
\def\moncompteur{\thesubsection.\arabic{prop}}
\def \texte #1 {\addtocounter{prop}1\noindent(\moncompteur) {\it #1}\\}
\newtheorem{theo}[prop]{Th{\'e}or{\`e}me}
\newtheorem{coro}[prop]{Corollaire}

\makeatletter

\@addtoreset{equation}{section}
\makeatother

\def \limind#1{\lim\limits_{\displaystyle\longrightarrow\atop {#1}}}
\def \limproj#1{\lim\limits_{\displaystyle\longleftarrow\atop {#1}}}
\def \limpro{\lim\limits_{\displaystyle\longleftarrow}}

\def \spec#1{\hbox{\rm Spec}\,(#1)}
\def \spf#1{\hbox{\rm Spf}\,(#1)}
\def \spm#1{\hbox{\rm Spm}\,(#1)}

\def \ker{\mathop{\hbox{\rm Ker}\,}}
\def \TS{\mathop{\hbox{\rm TS}\,}}

\def \Rlimproj#1{\RM\hspace{-0.1cm}\limproj{#1}}

\def \u#1{\underline{#1}}
\def \ov#1{\overline{#1}}
\def \surfleche#1{\buildrel#1\over\longrightarrow}
\def \dlog#1{\frac{d#1}{#1}}
\def \gsurfleche#1{\buildrel#1\over\longleftarrow}

\title{Classes de Chern en cohomologie rigide}
\author{Denis Petrequin }

\begin{document}
\maketitle
\tableofcontents
\bibliographystyle{plain}

\section*{Introduction}

La cohomologie rigide, introduite par Berthelot \cite{Berth8} est une th\'eorie cohomologique $p$-adique. C'est une g\'en\'eralisation de la cohomologie de Monsky-Washnitzer \cite{Mon-Was} et de la cohomologie cristalline \cite{Berth1}. Berthelot a d\'emontr\'e qu'elle v\'erifiait la propri\'et\'e de finitude \cite{Berth3}, la dualit\'e de Poincar\'e et la formule de K\"{u}nneth \cite{Berth7}. Dans cet article nous \'etudions les classes de Chern. Dans \cite{Pet2} nous nous interesserons au cas des classes de cycles. Nous obtenons que la cohomologie rigide satisfait tous les axiomes des cohomologies de Weil.

Pr\'ecisons les diff\'erentes parties de l'article.

Nous commen\c{c}erons par un chapitre pr\'eliminaire, dans lequel on trouvera des rappels sur la d\'efinition de la cohomologie rigide ainsi que sur ses propri\'et\'ees. 

\vspace{.3cm}

Dans le deuxi\`eme  chapitre nous construisons les classes de Chern pour
les vari{\'e}t{\'e}s propres. Du fait de la d{\'e}finition m{\^e}me de la
cohomologie rigide, cette construction s'apparente {\`a} une variante
{\it rigide} des constructions de classes de Chern en cohomologie de
De Rham. Cette m{\'e}thode a l'avantage d'{\^e}tre {\'e}l{\'e}mentaire et
explicite. Cependant, il semble difficile - ou au moins laborieux - de
d{\'e}montrer l'additivit{\'e} de classes de Chern ainsi construites.

\vspace{.3cm}

Dans le troisi{\`e}me chapitre, apr\`es avoir fait des rappels sur la cohomologie cristalline de niveau $m$, nous construisons des classes de Chern \`a valeur dans cette derni\`ere en nous basant sur la construction des classes de Chern en cohomologie cristalline de Berthelot et Illusie \cite{Berth-Ill}. Cela nous permet de r\'einterpreter la construction des classes de Chern rigides. Cette m\'ethode \`a l'avantage de donner une construction plus intrins\`eque ce qui nous permet d'utiliser les m\'ethodes classiques (restriction \`a un cas universel) pour d\'emontrer l'additivit\'e des classes de Chern. 

On donnera aussi des th\'eor\`emes de comparaison entre nos classes de Chern et les classes de Chern cristallines ou les classes de Chern \`a valeur dans le topos convergent d'Ogus \cite{Ogu} construites par Niziol \cite{Niz}. 

\vspace{.3cm}

Dans le quatri{\`e}me chapitre, nous nous int{\'e}ressons aux vari{\'e}t{\'e}s
non n{\'e}cessairement propres. Nous {\'e}tablissons quelques lemmes
relatifs aux faisceaux localement libres et aux compactifications. L'outil principal est le th\'eor\`eme de platification par \'eclatement de Raynaud-Gruson \cite{Gru-Ray}. Nous
d{\'e}montrons, par la suite, le th{\'e}or{\`e}me principal qui stipule que pour calculer
la premi{\`e}re classe de Chern d'un faisceau inversible $\LCC$ sur une
vari{\'e}t{\'e} $X$, il suffit de
trouver un faisceau inversible $\ov{\LCC}$ sur un compactifi{\'e} $\ov{X}$ de $X$ se
restreignant {\`a} $\LCC$ sur cette derni{\`e}re et de prendre
l'image de la premi{\`e}re classe de Chern de $\ov{\LCC}$ par le
morphisme de fonctorialit{\'e} $H^2_{rig}(\ov{X}) \to H^2_{rig}(X).$ Nous montrons en effet (th\'eor\`eme \ref{theo3}) que la classe ainsi trouv{\'e}e ne d{\'e}pend pas des choix faits.
Nous construisons ensuite les autres classes de Chern de mani{\`e}re
classique {\`a} l'aide de fibr{\'e}s projectifs.

\vspace{.3cm}

Cet article est la premi\`ere partie - l\'eg\`erement modifi\'ee - de ma th\`ese de doctorat \cite{moi}. Je tiens \`a remercier P. Berthelot pour tous les conseils qu'il m'a prodigu\'es tout au long de la r\'edaction de cet article.

\vspace{.5cm}

\section{Pr\'eliminaires et notations}
Tout au long de cet article, $k$ d\'esignera un corps de caract\'eristique $p > 0$. On appelera $k$-vari\'et\'e un sch\'ema s\'epar\'e de type fini sur $\spec{k}$.
On se donne un anneau de valuation discr{\`e}te complet
d'in{\'e}gale caract{\'e}ristique $\VCC$ de corps r{\'e}siduel $k$. On note
alors $K$ son corps des fractions.

Rappelons pour commencer la construction de la cohomologie rigide \cite{Berth3}. Soit $X$ une $k$-vari\'et\'e. Il existe d'apr\`es Nagata \cite{Nag}, une vari\'et\'e propre $\ov{X}$ et une immersion ouverte $j: X \hookrightarrow \ov{X}$.  On suppose alors qu'il existe une immersion ferm{\'e}e $\ov{X} \hookrightarrow \YCC$ dans un sch{\'e}ma formel $\YCC$ sur $\spf{\VCC}$ lisse au voisinage de
$X$ - cette condition technique peut \^etre supprim\'ee en utilisant des r\'esolutions de Cech, nous la garderons pour simplifier le propos. On consid{\`e}re la fibre g{\'e}n{\'e}rique rigide $\YS$ de $\YCC$ \cite{Ray} et on note $sp : \YS \to \YCC$ le morphisme de sp\'ecialisation. Rappelons \cite[1.1]{Berth5} que si on note $\YCC_0$ la r\'eduction de $\YCC$ sur $k$, on appelle, pour tout sous-$k$-sch\'ema $T$ de $\YCC_0$, {\it tube de $T$} et on note $]T[$ le sous-ensemble $sp^{-1}(T)$ des points de $\YS$ qui se sp\'ecialisent dans $T$.

Avec les notations pr\'ec\'edentes, on appelle voisinage strict \cite[1.2]{Berth5} de $]X[$ dans $]\ov{X}[$, tout ouvert $V$ de $]\ov{X}[$ tel que le recouvrement $(V, ]\ov{X} - X[)$ soit admissible. D\`es lors, pour tout faisceau $\ECC$ sur $]\ov{X}[$, on note
$$j^\dag \ECC := \limind{V}\alpha_{V*}\alpha_V^* \ECC,$$
o\`u la limite inductive est prise sur tous les voisinages stricts de $]X[$ dans $]\ov{X}[$ et $\alpha_V$ d\'esigne l'inclusion $V \hookrightarrow ]\ov{X}[$.

On sait alors \cite{Berth5} que le complexe $\RM sp_* j^\dag \Omega^\star_{]\ov{X}[}$ vu comme objet de la cat\'egorie d\'eriv\'ee des complexes de $K$-vectoriels sur $\ov{X}$ est ind\'ependant \`a isomorphisme canonique pr\`es du choix de $\YCC$. On posera donc :

$$\RM\u{\Gamma}_{rig}((X,\ov{X})/K) := \RM sp_* j^\dag \Omega^\star_{]\ov{X}[}.$$

Berthelot montre alors que les groupes de cohomologie 
$$H^i(\ov{X}, \RM sp_* j^\dag \Omega^\star_{]\ov{X}[})$$
ne d\'ependent pas de la compactification $\ov{X}$ choisie.

On pose donc

$$H^i_{rig}(X/K) := H^i(\ov{X}, \RM sp_* j^\dag \Omega^\star_{]\ov{X}[}).$$

Quand il n'y aura pas d'ambigu\"it{\'e} sur la compactification nous
noterons $\RM\u{\Gamma}_{rig}(X/K)$ pour $\RM sp_* j^\dag \Omega^\star_{]\ov{X}[}.$

\section{Cas des vari\'et\'es propres}
\label{propre}
On va se restreindre au cas des vari\'et\'es propres. Nous allons construire la premi\`ere classe de Chern d'un faisceau inversible \`a l'aide d'un calcul de cocycle.

\rqe Pour une $k$-vari\'et\'e propre, la cohomologie rigide est isomorphe \`a la cohomologie du topos convergent d'Ogus \cite{Ogu}. Nos calculs ne sont donc qu'une r\'einterpr\'etation de la construction des classes de Chern \`a valeurs dans le topos convergent \cite[A]{Niz} (proposition \ref{conv}). Cependant ce calcul de cocycle sera n\'ecessaire pour pouvoir \'etudier par la suite le cas des vari\'et\'es ouvertes.

\subsection{Premi\`ere classe de Chern}
\label{unrefdeplus}
Soient $X$ une vari\'et\'e propre, $\FCC$ un faisceau inversible sur $X$ et une immersion
ferm{\'e}e $X \hookrightarrow \YCC$ dans un sch\'ema formel sur $\spf{\VCC}$ lisse au voisinage de $X$.
\begin{lemme}
\label{debile}
Il existe un recouvrement affine $\UG=(\UCC_i)_{i\in \Lambda}$ de $\YCC$ tel
que si on note $\UG_X$ le recouvrement induit par $\UG$ sur $X$, le
faisceau $\FCC$ soit trivialis{\'e} sur $\UG_X$.
\end{lemme}
\dem Quitte {\`a} prendre un recouvrement affine de $\YCC$ on peut
supposer que $\YCC$ et $X$ sont affines. On pose alors $\YCC =
\spf{\ACC}$ et $X =\spec{A}$ o{\`u} $A=\ACC/I.$ On choisit un
recouvrement $\UG = (X_i)_{i\in \Lambda}$ de $X$ trivialisant $\FCC$. Quitte {\`a}
raffiner le recouvrement, on peut supposer que pour tout $i \in \Lambda$, il
existe $f_i \in A$ tel que $X_i = D(f_i)$. Pour tout $i \in \Lambda$, on
choisit  alors $\widetilde{f_i}$ un
rel{\`e}vement de $f_i$ dans $\ACC$, les
ouverts $\spf{\widehat{\ACC_{\widetilde{f_i}}}}$ recouvrent alors l'image de $X$ dans $\YCC$. En rajoutant des ouverts de $\YCC - X$,
on obtient le recouvrement voulu.

\findem

On choisit donc un tel recouvrement $\UG =(\UCC_i)_{i\in \Lambda}$. D{\`e}s lors
on
consid{\`e}re $\UG_X = (X_i)_{i\in \Lambda}$ le recouvrement induit sur
$X$. Pour tout $i \in I$, on notera $\UCC_i = \spf{\ACC_i}$, $X_i =
\spec{A_i}$ o{\`u} $A_i = \ACC_i/I_i$. 
On se donne alors un cocycle $(u) \in Z^1(\UG_X, \OCC_X^*)$
repr{\'e}sentant la classe du faisceau $\FCC$ dans $H^1(\UG_X, \OCC_X^*).$

Commen\c{c}ons par un lemme :

\begin{lemme}
\label{laumon}
Soit $\YCC = \spf{\ACC}$ un $\VCC$-sch\'ema formel affine. Soit $X=\spec{A}$ un ferm\'e de la fibre sp\'eciale $\YCC_0$ de $\YCC_0$. On note $]X[$ le tube de $X$ dans la fibre g\'en\'erique rigide $\YS=\YCC_K$. Pour tout $\widetilde{u} \in \ACC$ relevant un \'el\'ement inversible $u$  de $A$, la restriction de $u$ \`a $]X[$ est inversible.
\end{lemme}

\dem On note $v$ l'inverse de $u$ dans $A$ et on choisit $\widetilde{v}$ un rel\`evement de $v$ dans $\ACC$. D\`es lors, si on note $I$ le noyau de l'application $\ACC \to A$, il existe $a \in I$ tel que 
$$\widetilde{u}.\widetilde{v} = 1+a.$$
On regarde $a$ comme une fonction analytique sur $\YS$. Par d\'efinition du tube $]X[$, on a que, pour tout $x \in ]X[$ :
$$|a(x)| <1.$$
La fonction $1+a$ est donc inversible sur cet ouvert et son
inverse est la s{\'e}rie :
$$\sum_{n=0}^\infty (-1)^na(x)^n$$
qui converge.

L'inverse de $\widetilde{u}$ est alors 
$$\widetilde{u}^{-1} = \widetilde{v}(1+a)^{-1}.$$
\findem

Pour tous $i,j$ des {\'e}l{\'e}ments de $\Lambda$ on note $v_{ij} \in A_{ij}$
l'inverse de $u_{ij}$. On choisit alors des rel{\`e}vements de ces deux
{\'e}l{\'e}ments que l'on note ; 
$$\widetilde{u}_{ij} \in \ACC_{ij} \hbox{ et } \widetilde{v}_{ij} \in \ACC_{ij}.$$

On sait gr\^ace au lemme \ref{laumon} que 
$$\widetilde{u}_{ij} \in \Gamma(]X_{ij}[ ,\OCC_{]X[}^*).$$

De plus pour tous $i, j, k$ {\'e}l{\'e}ments de $\Lambda$, on regarde :
$$\widetilde{u}_{ij}.\widetilde{v}_{ik}.\widetilde{u}_{jk}.$$

Comme $(u)$ est un cocycle, la r{\'e}duction modulo $I_{ijk}$ de ce terme est
{\'e}gale $1$. Il existe donc $b_{ijk} \in I_{ijk}$ tel que :
$$\widetilde{u}_{ij}.\widetilde{v}_{ik}.\widetilde{u}_{jk}=1+b_{ijk}.$$
A partir de l\`a, dans $\Gamma(]X_{ijk}[, \OCC_{]X[})$, il existe $y_{ijk}$ tel que
$$\widetilde{u}_{ij}.\widetilde{u}_{ik}^{-1}.\widetilde{u}_{jk} = (1+b_{ijk}).(1+a_{ij})^{-1} = 1+y_{ijk}.$$
De plus, il est clair que pour tout $x \in ]X_{ijk}[$, on a :
$$|y_{ijk}(x)| < 1.$$ 

Comme ci dessus :
$$c_{ijk}:=\log(\widetilde{u}_{ij}.\widetilde{u}_{ik}^{-1}.\widetilde{u}_{jk})=\sum_{n=1}^\infty
(-1)^{n-1}\frac{y_{ijk}^n}{n}$$
converge sur $]X_{ijk}[$.

On note $\UG_K$ le
recouvrement admissible de $]X[$ par les ouverts $]X_{i}[$ \cite[1.1.14]{Berth5}. 
On consid{\`e}re alors le bicomplexe $C^\star(\UG_K,
\Omega^\star_{]X[})$. On note $d$ la diff{\'e}rentielle provenant de
la diff{\'e}rentielle du complexe $\Omega^\star_{]X[}$ et $\delta$
celle qui provient du complexe de Cech. On consid{\`e}re alors le
complexe simple associ{\'e}, encore not{\'e} $C^\star(\UG_K,
\Omega^\star_{]X[})$, muni de sa diff{\'e}rentielle d{\'e}finie sur le
terme $C^p(\UG_K, \Omega_{]X[}^q)$ par :
$$\Delta = \delta + (-1)^pd.$$

On construit un {\'e}l{\'e}ment 
$$c_1(u) \in C^2(\UG_K,
\Omega^\star_{]X[}) = C^0(\UG_K,\Omega^2_{]X[}) \oplus
C^1(\UG_K,\Omega^1_{]X[}) \oplus C^2(\UG_K,\OCC_{]X[})$$ 
en posant :
$$c_1(u)_i = 0,$$
$$c_1(u)_{ij} := d(\widetilde{u}_{ij}).\widetilde{u}_{ij}^{-1} \in C^1(\UG_K,
\Omega^1_{]X[}),$$
et
$$c_1(u)_{ijk} :=
-\log(\widetilde{u}_{ij}.\widetilde{u}_{ik}^{-1}.\widetilde{u}_{jk}) \in C^2(\UG_K, \OCC_{]X[}).$$

Il est clair qu'avec les notations ci-dessus on a :
$$\Delta(c_1(u))=0.$$

Notons que l'{\'e}l{\'e}ment ainsi construit d{\'e}pend du choix de nos
rel{\`e}vements ; nous allons montrer cependant que sa classe dans
$\HM^2(\UG_K,\Omega^\star_{]X[})$ n'en d{\'e}pend pas.

Soient  $(u) \in C^1(\UG_X, \OCC_X^*)$  et  $(\theta) \in C^0(\UG_X, \OCC_X^*)$. On pose :
\begin{equation}
\label{truc}
u'=u.\delta(\theta).
\end{equation}
Comme pr{\'e}c{\'e}demment, on choisit pour tout $i \in \Lambda$, des
rel{\`e}vements $\widetilde{\theta}_i \in \ACC_i$ de $\theta_i$. Le lemme \ref{laumon} nous dit que la fonction $\widetilde{\theta}_i$
est inversible  sur $]X_i[$. On pose :
\begin{equation}
\label{zeta1}
\zeta_i = \frac{d\widetilde{\theta}_i}{\widetilde{\theta}_i} \in C^0(\UG_K,
\Omega^1_{]X[}).
\end{equation}
De plus l'{\'e}quation \ref{truc} nous dit que pour tous $i,j \in \Lambda$, si
$\widetilde{u}_{ij}$ et $\widetilde{u}'_{ij}$ sont respectivement des
rel{\`e}vements de $u_{ij}$ et $u'_{ij}$, il
existe $\alpha_{ij} \in \Gamma(]X_{ij}[, \OCC_{]X[})$  tel que pour tout $x \in ]X_{ij}[$ on ait $|\alpha_{ij}(x)| <1$ et que :
$$\widetilde{\theta}_i \widetilde{u}'_{ij} = \widetilde{\theta}_j \widetilde{u}_{ij}(1+\alpha_{ij})$$
o{\`u} $\widetilde{\theta}_i$ et $\widetilde{\theta}_j$ sont vus comme des
{\'e}l{\'e}ments de $\Gamma(]X_{ij}[, \OCC_{]X[})$ par les fl{\`e}ches de restriction {\'e}videntes.

On pose alors :
\begin{equation}
\label{zeta2}
\zeta_{ij} = -\log(1+\alpha_{ij}) \in C^1(\UG_K, \OCC_{]X[}),
\end{equation}
o\`u le $\log$ est d\'efini comme pr\'ec\'edemment.

Les formules \ref{zeta1} et \ref{zeta2} d{\'e}finissent :
$$\zeta \in C^1(\UG_K,\Omega^\star_{]X[}) =
C^0(\UG_K,\Omega^1_{]X[}) \oplus C^1(\UG_K,\OCC_{]X[}).$$
\begin{lemme}
\label{bordel}
 Avec les notations ci-dessus, on a :
$$c_1(u')-c_1(u) = \Delta(\zeta).$$
\end{lemme}

\dem
On a le diagramme suivant
$$\xymatrix{ & C^0(\UG_K,\Omega^1_{]X[}) \ar[dl]_{d} \ar[dr]^{\delta}
  & &
  C^1(\UG_K,\OCC_{]X[}) \ar[dl]_{-d} \ar[dr]^{\delta} & \\
C^0(\UG_K,\Omega^2_{]X[}) &  & C^1(\UG_K,\Omega^1_{]X[}) & & C^2(\UG_K,\OCC_{]X[}).}$$
On va calculer $\Delta(\zeta).$ 
On a :
$$\Delta(\zeta)_i =
d\left(\frac{d\widetilde{\theta}_i}{\widetilde{\theta}_i}\right) =
0.$$
$$\begin{array}{ccc} \Delta(\zeta)_{ij} & = & \delta(\zeta_{i}) -
  d(\zeta_{ij}) \\
& = &
  \frac{d\widetilde{\theta}_j}{\widetilde{\theta}_j}-\frac{d\widetilde{\theta}_i}{\widetilde{\theta}_i} + d\log(1+\alpha_{ij}) \\
  & = & \frac{d\widetilde{u'}_{ij}}{\widetilde{u'}_{ij}}
  - \frac{d\widetilde{u}_{ij}}{\widetilde{u}_{ij}}. \end{array}$$
$$\begin{array}{ccc} \Delta(\zeta)_{ijk} & = & \delta(\zeta_{ij}) \\
& = &
-\log\left(\frac{\widetilde{\theta}_i\widetilde{u'}_{ij}.\widetilde{\theta}_k\widetilde{u}_{ik}.\widetilde{\theta}_j\widetilde{u'}_{jk}}{\widetilde{\theta}_j\widetilde{u}_{ij}.\widetilde{\theta}_i\widetilde{u'}_{ik}.\widetilde{\theta}_k\widetilde{u}_{jk}}\right)
\\
& = &
\log(\widetilde{u}_{ij}.\widetilde{u}_{ik}^{-1}.\widetilde{u}_{jk}) - \log(\widetilde{u'}_{ij}.\widetilde{u'}_{ik}^{-1}.\widetilde{u'}_{jk}) .
\end{array}$$
\findem
On a donc montr\'e que l'\'el\'ement $c_1(u) \in \HM^2(\UG_K,\Omega^\star_{]X[})$ ne d\'epend ni du repr\'esentant du cocycle choisi ni des choix de rel\`evements (ce dernier cas est obtenu en prenant pour $\theta$ le cocycle trivial : $\theta_i = 1$).

On construit donc bien ainsi une application
$$H^1(X,\OCC_X^*) \to \HM^2(\UG_K,\Omega^\star_{]X[}).$$

Il est clair que cette application est compatible au raffinement du recouvrement $\UG$.

De plus, d'apr\`es le th\'eor\`eme $B$ \cite{Kie}, on sait que :
$$H^2_{rig}(X/K)=\HM^2(]X[,\Omega^\star_{]X[})=\HM^2(\UG,\Omega^\star_{]X[}).$$

Gr\^ace a un plongement diagonal, on peut alors montrer que l'application ainsi construite ne d\'epend pas du choix du plongement.
On note donc :
$$c_{1,rig} : H^1(X,\OCC_X^*) \to H^2_{rig}(X/K),$$
l'application obtenue.

\begin{prop}[Fonctorialit{\'e}]
Soient $X$ et $X'$ deux $k$-vari{\'e}t{\'e}s propres et $f : X' \to X$
un morphisme. Pour tout faisceau inversible $\FCC$ sur $X$ on a :
$$c_{1,rig}(f^*\FCC) = f^*c_{1,rig}(\FCC).$$
\end{prop}
\dem La d{\'e}monstration se fait par un calcul classique. Nous n'allons 
cependant pas le d\'evelopper ici car nous allons donner plus tard
une d{\'e}finition cristalline des classes de Chern et la
fonctorialit{\'e} sera alors {\'e}vidente.

\begin{prop}[Multiplicativit\'e]
\label{multi}
Soient $X$ une $k$-vari\'et\'e propre et deux faisceaux inversibles $\FCC$ et $\FCC'$ sur $X$, on a :
$$c_{1,rig}(\FCC \otimes \FCC') = c_{1,rig}(\FCC) + c_{1,rig}(\FCC').$$
L'application $c_{1,rig}$ est donc un morphisme de groupe.

\end{prop}

\dem On choisit un recouvrement $\UG$ qui trivialise les deux faisceaux. D\`es lors en notant $(u)$ et $(u')$ des cocycles repr\'esentant respectivement $\FCC$ et $\FCC'$, le faisceau $\FCC \otimes \FCC'$ est repr\'esent\'e par le cocycle $(u.u')$.
En choisissant $(\widetilde{u})$ et $(\widetilde{u}')$ des rel\`evements de $(u)$ et $(u')$ respectivement, on a :
$$c_1(u.u')_{ij} = d(\widetilde{u}_{ij}.\widetilde{u}'_{ij}).(\widetilde{u}_{ij}.\widetilde{u}'_{ij})^{-1} = c_1(u)_{ij} + c_1(u')_{ij}.$$
et
$$c_1(u.u')_{ijk} = - \log(\widetilde{u}_{ij}\widetilde{u'}_{ij}.\widetilde{u}_{ik}^{-1}\widetilde{u'}_{ik}^{-1}.\widetilde{u}_{jk}\widetilde{u}'_{jk}) = c_1(u)_{ijk}+c_1(u')_{ijk}.$$

\findem

\begin{prop}[Extension des scalaires]
Soit $K'$ une extension de $K$, d'anneau des entiers $\VCC'$ et de corps r\'esiduels $k'$. Pour tout faisceau inversible $\FCC$ sur $X$ on note $\FCC'$ le faisceau inversible sur $X' = X\times_k k'$ obtenu par changement de base. Le morphisme
$$H^2_{rig}(X/K) \to H^2_{rig}(X'/K')$$
envoie alors $c_{1,rig}(\LCC)$ sur $c_{1,rig}(\LCC')$.

\end{prop}

\dem Il est clair que toutes nos constructions commutent aux changements de corps de base.
\findem
Par la suite, on pourra donc choisir pour $\VCC$ un anneau de Cohen associ\'e \`a $k$.

\subsection{Les autres classes}
Nous allons suivre la m{\'e}thode classique \cite{Gro3} afin de
construire les autres classes de Chern d'un faisceau localement libre.

Nous allons calculer la cohomologie rigide d'un fibr{\'e} projectif sur
une base propre. Plus pr{\'e}cisement on se donne $X$ une
$k$-vari{\'e}t{\'e} propre et $\ECC$ un faisceau localement libre de rang $r$. On note alors
$\PM=\PM(\ECC)$ le fibr{\'e} projectif associ{\'e}, $p : \PM \to X$ la
projection et 
$$\xi=c_{1,rig}(\OCC_{\PM}(1)) \in H_{rig}^2(\PM/K).$$
On peut alors regarder $\xi$ comme un morphisme dans $D^+(X_{Zar})$ :
$$\xi : \ZM_X \to \RM p_* \RM\u{\Gamma}_{rig}(\PM/K)[2],$$
o\`u $\ZM_X$ est le faisceau constant $\ZM$ sur $X$.
Par cup-produit on a  alors :
$$\xi^i : \ZM_X \to \RM p_* \RM\u{\Gamma}_{rig}(\PM/K)[2i].$$
Enfin on a le morphisme de fonctorialit{\'e} :
$$\RM\u{\Gamma}_{rig}(X/K) \to \RM p_*\RM\u{\Gamma}_{rig}(\PM/K).$$
On en d{\'e}duit alors le morphisme dans $D^+(X_{Zar})$ :
\begin{equation}
\label{bigoplus}
\bigoplus_{i=0}^{r-1} \xi^i : \bigoplus_{i=0}^{r-1}\RM\u{\Gamma}_{rig}(X/K)[-2i]
\to \RM p_* \RM\u{\Gamma}_{rig}(\PM/K)
\end{equation}

\begin{prop}
Avec les notations pr{\'e}c{\'e}dentes, l'application (\ref{bigoplus}) est
un isomorphisme.
\label{struc}
\end{prop}

\dem La d{\'e}monstration est classique. La question
{\'e}tant locale, on peut supposer que $\ECC = \OCC_X^r$. On se donne
alors une immersion ferm{\'e}e $X \hookrightarrow \YCC$ dans un
$\VCC$-sch{\'e}ma formel lisse. On a alors le diagramme commutatif suivant :
$$\xymatrix{\PM^r_X \ar[r] \ar[d] & \PM^r_{\YCC} \ar[d] & \PM^r_{]X[}
  \ar[d] \ar[l] \\
X \ar[d] \ar[r] & \YCC \ar[d]  & ]X[ \ar[d] \ar[l]\\ \spec{k} \ar[r] & \spf{\VCC} & \spm{K} \ar[l] }$$
Le faisceau $\OCC_{\PM^r_X}(1)$ se relevant en
$\OCC_{\PM^r_{\YCC}}(1)$, on est ramen{\'e} au lemme suivant :

\begin{lemme}
\label{18}
Soit $X$ un affinoide, la fl\`eche d\'efinie comme pr\'ec\'edemment :
$$\bigoplus_{i=0}^{r-1} \Omega^\star_{X/K}[-2i] \to \RM p_* \Omega^\star_{\PM_X/K}$$
est un isomorphisme
\end{lemme}

\dem La d\'emonstration est analogue au cas analytique complexe \cite{Verd}. Ce r\'esultat peut \^etre obtenu directement en utilisant un th\'eor\`eme de type GAGA relatif en g\'eom\'etrie rigide \cite[2.8]{Lutk}. Il suffit pour cela de voir que la d\'emonstration de \cite{Ser} se recopie dans notre cas.

\findem

\begin{coro} 

\label{juste}
Avec les notations ci-dessus, on a pour tout $n$ la
d{\'e}composition suivante :
$$H^n_{rig}(\PM/K)=\bigoplus_{i=0}^{r-1}H^{n-2i}_{rig}(X/K).\xi^i.$$

\end{coro}

On d{\'e}finit les classes de Chern sup{\'e}rieures comme dans le cas
classique. On applique la d{\'e}composition du corollaire \ref{juste} {\`a}
$\xi^r$ et on obtient :
\begin{equation}
\label{decompo}
\xi^r=\sum_{i=1}^r (-1)^{i+1} c_i(\ECC)\xi^{r-i},
\end{equation}
avec
$$c_i(\ECC) \in H^{2i}_{rig}(X/K).$$

\begin{Defi}
Avec les notations pr{\'e}c{\'e}dentes, pour $0 < i \leq r$ on appelle $i$-{\`e}me classe de
Chern de $\ECC$ la classe $c_i(\ECC)$. On pose de plus 
$$c_0(\ECC) = 1 \hbox{ dans } H^0_{rig}(X/K).$$
\end{Defi} 

\rqe En appliquant ce que l'on vient de voir \`a un fibr\'e inversible $\FCC$, on retrouve notre premi\`ere classe de Chern. En effet la d\'ecomposition \ref{decompo} devient alors :
$$c_{1, rig}(\OCC_{\PM(\FCC)}(1)) = p^*c_1(\FCC).$$
Or la fonctorialit\'e de $c_{1, rig}$ nous donne
$$c_{1, rig}(\OCC_{\PM(\FCC)}(1)) = p^*c_{1, rig}(\FCC).$$
On conclut en utilisant que $p^*$ est injectif.

\begin{prop}[Fonctorialit\'e]
Soient $X$ et $X'$ deux $k$-vari{\'e}t{\'e}s propres et $f : X' \to X$
un morphisme. Pour tout $\ECC$ faisceau localement libre de rang $r$ sur $X$ et tout $i$, on a 
$$c_{i}(f^*\ECC) = f^*c_{i}(\ECC).$$
\end{prop}

\dem
On note $\PM' = \PM(f^*(\ECC))$ et $\xi' = c_1(\OCC_{\PM'}(1))$. La fonctorialit\'e de la premi\`ere classe de Chern nous assure que 
$$\xi' = f^* \xi.$$
 On applique alors $f^*$ \`a la d\'ecomposition \ref{decompo}. Le morphisme $f^*$ \'etant compatible au produit, on obtient dans $H^{2r}_{rig}(\PM'/K)$ :
$${\xi'}^r = \sum_{i=1}^r (-1)^i f^*c_i(\ECC){\xi'}^{r-i}.$$
La d\'ecomposition \'etant unique on a :
$$c_i(f^*\ECC) = f^*c_i(\ECC).$$

\findem

\rqe La d\'efinition utilis\'ee pour la cohomologie rigide ne permet pas d'utiliser la m\'ethode classique \cite{Gill} pour d\'emontrer l'additivit\'e des classes de Chern. Nous diff\'erons donc la d\'emonstration de cette propri\'et\'e \`a la section suivante.

\section{Construction cristalline et comparaisons}

Apr\`es avoir fait quelques rappels sur les topos cristallins de niveau $m$, nous construisons des classes de Chern \`a valeurs dans ces derniers en nous inspirant du cas cristallin classique (de niveau $0$) \cite{Berth-Ill}. Cela nous permettera de donner une construction cristalline des classes de Chern rigides construites pr\'ec\'edemment. Cette formulation permet d'appliquer la m\'ethode classique \cite{Gill, Berth-Ill} pour d\'emontrer l'additivit\'e des classes de Chern en se ramenant \`a une situation universelle.

On \'etablira aussi un th\'eor\`eme de comparaison avec les classes de Chern \`a valeur dans le topos convergent \cite{Niz}. 

\subsection{Les topos $m$-cristallins}
Nous allons {\'e}noncer les principaux r{\'e}sultats sur la cohomologie
cristalline de niveau $m$. Pour de plus amples informations on pourra
consulter \cite{Tri}. Il sera aussi utile pour les puissances divis{\'e}es
de niveau $m$ de regarder \cite{Berth4}.

On fixe une fois pour toutes un nombre premier $p$. Commen\c{c}ons par rappeler la d\'efinition d'un $m$-PD-id\'eal et de l'enveloppe \`a puissances divis\'ees partielles de niveau $m$. On trouvera un expos\'e d\'etail\'e de cette th\'eorie dans \cite[1]{Berth4}.

\begin{Defi}
Soient $A$ une $\ZM_{(p)}$-alg\`ebre, $I \subset A$ un id\'eal et $m \geq 0$ un entier. On appelle structure partielle d'id\'eal \`a puissances divis\'ees de niveau $m$ sur $I$ ($m$-PD-structure) la donn\'ee d'un PD-id\'eal $(J,\gamma) \subset I$ tel que
$$I^{(p^m)} +pI \subset J.$$
\end{Defi} 

\begin{Defi}
Soient $A$ une $\ZM_{(p)}$-alg\`ebre et $I \subset A$ un id\'eal quelconque. Il existe une $A$-alg\`ebre $P^m(I)$ et un $m$-PD-id\'eal $\ov{I} \subset P^m(I)$ tel que $IP^m(I) \subset \ov{I}$, qui soient universels pour les morphismes de $A$ dans un anneau $A'$ envoyant $I$ dans un $m$-PD-id\'eal $I'$. Cette alg\`ebre munie de son $m$-PD-id\'eal est appel\'ee l'enveloppe \`a puissances divis\'ees partielles de niveau $m$ de $(A,I)$.
\end{Defi}

\rqe On a le diagramme commutatif suivant :
$$\xymatrix{ & P^m(I) \ar[dr] & \\ A \ar[ur] \ar[rr] & & A/I.}$$

Soit $C$ un anneau de Cohen pour $k$. On rappelle \cite[I.1.2.4]{Berth1} qu'il existe une unique structure de $PD$-id\'eal sur l'id\'eal maximal de $C$ not\'ee $\gamma$. On note $S = \spf{C}$. Dans toute la suite nous ne regarderons que les topos cristallins de base $S$ car c'est eux qui interviennent dans la comparaison avec la cohomologie rigide. 

\begin{Defi}
Soient $X$ un $S$-sch{\'e}ma sur lequel $p$ est localement nilpotent et $m \in
\NM$. On appelle site cristallin de niveau $m$ de $X$ relativement {\`a}
$(S, (p), \gamma)$ et on note $Cris^m(X/S)$ le site ayant pour objets les
$S$-immersions ferm{\'e}es $U \hookrightarrow T$ o{\`u} $U$ est un ouvert de $X$,  $p$ est localement nilpotent sur $T$ et o{\`u} l'id{\'e}al d{\'e}finissant ces immersions est muni d'une $m$-PD-structure compatible {\`a} $\gamma$.
\end{Defi}

Par la suite nous appellerons {\it topos cristallin de niveau $m$} et nous
noterons $(X/S)_{cris}^m$ le topos associ{\'e} au site $Cris^m(X/S)$.

\rqe pour $m=0$ on retrouve le topos cristallin classique.

Il est possible de donner une description des objets du topos cristallin de niveau
$m$ similaire {\`a} celle que l'on a pour le topos cristallin
classique. Pour d{\'e}finir un faisceau $\FCC$ du topos cristallin de niveau $m$, il suffit de donner pour tout {\'e}l{\'e}ment $(U,T)$ du site $m$-cristallin, un faisceau zariskien $\FCC_{(U,T)}$ (not\'e aussi $\FCC_T$) sur $T$ ; ces donn{\'e}es {\'e}tant assujetties {\`a} des hypoth{\`e}ses de compatibilit{\'e}s similaires {\`a} celles du cas classique.

On peut de cette mani\`ere construire les faisceaux suivants :
\begin{itemize}
\item Le faisceau $\OCC_{X/S}^m$ est d{\'e}fini par :
$$(\OCC_{X/S}^m)_{(U,T)}=\OCC_T.$$
\item 
 Le faisceau $\ICC_{X/S}^m$ est d{\'e}fini par :
$$(\ICC_{X/S}^m)_{(U,T)}=\ICC_T$$
o{\`u} $\ICC_T$ est le faisceau d{\'e}finissant l'immersion $U
\hookrightarrow T$.
\item 
 Le faisceau $\JCC_{X/S}^m$ est d{\'e}fini par :
$$(\JCC_{X/S}^m)_{(U,T)}=\JCC_T$$
o{\`u} $\JCC_T$ est le PD-id{\'e}al d{\'e}finissant la $m$-PD structure sur $\ICC_T$.
\end{itemize}

\notation Pour tout $k$-vari\'et\'e $X$ et tout $m$, on note :
$$H^i_{m-cris}(X) := H^i((X/S)^m_{cris}, \OCC_{X/S}^m).$$

Nous allons maintenant {\'e}tudier comment se comportent ces topos quand
on fait varier $m$.

\begin{prop}
Pour tous $m, m'$ avec $m \leq m'$ il existe un morphisme de topos :
$$i_{m',m} : (X/S)_{cris}^m \to (X/S)_{cris}^{m'}.$$
De plus ces morphismes v{\'e}rifient que pour tout triplet $m \leq m'
\leq m''$ on a :
$$ i_{m'',m} \iso i_{m'',m'} \circ i_{m',m}.$$
\end{prop}

Ces morphismes sont donn{\'e}s de la mani{\`e}re suivante :
\begin{itemize}

\item pour tout $\ECC'$ dans $(X/S)_{cris}^{m'}$ et tout $(U,T,\ICC) \in
  Cris^{m}(X/S)$ on a :
$$\Gamma((U,T),{i_{m',m}^*}\ECC') = \Gamma((U,T),\ECC')$$
en utilisant que l'id\'eal $\JCC$ d\'efinissant la $m$-PD-structure de $\ICC$ d\'efinit aussi une $m'$-PD-structure sur $\ICC$.

\item pour tout $\ECC$ dans $(X/S)_{cris}^m$ et tout $(U,T,\ICC) \in
  Cris^{m'}(X/S)$ on a :
$$\Gamma((U,T),{i_{m',m}}_*\ECC) = \Gamma((U,T^m),\ECC)$$
o{\`u} $T^m$ est l'enveloppe {\`a} puissances divis{\'e}es de niveau $m$ de
$\ICC$ compatible aux puissances divis{\'e}es d{\'e}finissant la $m'$-PD-structure de $\ICC$.
\end{itemize}

On construit alors un morphisme canonique
$$u_{m',m} : \OCC_{X/S}^{m'} \to i_{m',m*}\OCC_{X/S}^m$$
qui est donn{\'e} sur un {\'e}paississement $(U,T)$ par le morphisme canonique $\OCC_T \to \OCC_{T^m}.$
Le morphisme $i_{m',m}$ devient ainsi un morphisme de topos annel\'es.

Pour finir avec ces g\'en\'eralit\'es sur les topos $m$-cristallins, notons qu'il existe, comme dans le cas classique, un morphisme de projection sur le topos zariskien :
$$u_{X/S}^m : (X/S)_{cris}^m \to X_{Zar},$$
ainsi qu'un morphisme d'inclusion du topos Zariskien dans le topos $m$-cristallin :
$$i_{X/S}^m : X_{Zar} \to (X/S)_{cris}^m.$$
\begin{lemme}
\label{commut}
Pour tous $m, m'$ avec $m \leq m'$ on a les diagrammes commutatifs suivants :
$$\xymatrix{(X/S)^m_{cris} \ar[d]_{i_{m',m}} \ar[r]^{u_{X/S}^m} & X_{Zar} \ar@{=}[d] \\ (X/S)^{m'}_{cris}  \ar[r]^{u_{X/S}^{m'}} & X_{Zar}. } \hbox{\hspace{.5cm} et \hspace{.5cm}}\xymatrix{X_{Zar} \ar@{=}[d] \ar[r]^{i_{X/S}^m} & (X/S)^m_{cris} \ar[d]_{i_{m',m}}  \\ X_{Zar} \ar[r]^{i_{X/S}^{m'}} &   (X/S)^{m'}_{cris} . } $$
\end{lemme}

\subsection{Classes de Chern $m$-cristallines}
Nous allons g\'en\'eraliser la construction de \cite{Berth-Ill} au cas des topos $m$-cristallins.
On se fixe $m \in \NM$, avec les notations pr\'ec\'edentes, on a ,comme dans le cas de la cohomologie cristalline de niveau $0$, une
suite exacte dans $(X/S)_{cris}^m$ : 
$$0 \to 1 + \ICC_{X/S}^m \to (\OCC_{X/S}^m)^* \to
i_*^m(\OCC_{X}^*) \to 0.$$
On consid{\`e}re le cobord de cette suite exacte qui nous fournit un
morphisme dans $D^+((X/S)^m_{cris},Ab)$ :

\begin{equation}
\label{m-prem}
i^m_*(\OCC_{X}^*)[1] \to (1+\ICC_{X/S}^{m})[2].
\end{equation}

L'{\'e}l{\'e}vation  {\`a} la puissance $p^{m}$ nous
donne un morphisme :
$$1 +\ICC^m_{X/S} \to 1 + \JCC^m_{X/S}.$$

L'id{\'e}al $\JCC^m_{X/S}$ est muni d'une structure de
$PD$-id{\'e}al, 
on peut donc construire un logarithme. En le composant avec le
morphisme pr{\'e}c{\'e}dent on obtient :
$$\begin{array}{cccc} \psi_m : & 1 + \ICC^m_{X/S}
  & \to &  \JCC^m_{X/S}\\
    & (1+x) & \mapsto & \log((1+x)^{p^{m}})
\end{array}$$

En composant \ref{m-prem} avec $\psi_m$, on obtient

\begin{equation}
\label{m-c1}
c_{1,m} : i^m_*(\OCC_{X}^*)[1] \to \JCC^m_{X/S}[2].
\end{equation}

Par suite en composant avec $\JCC^m_{X/S} \to \OCC_{X/S}^m$ et en passant \`a la cohomologie on a :
\begin{equation}
\label{m-c1g}
c_{1,m} : H^1(X,\OCC_X^*) \to H^2_{m-cris}(X).
\end{equation}

\begin{prop}[Fonctorialit{\'e}]
\label{m-fonctorielle}
Soient $X$ et $Y$ deux $k$-vari{\'e}t{\'e}s, $f : X \to
Y$ un morphisme et $\FCC$ un faisceau inversible sur $Y$, on
a
$$c_1(f^*\FCC)=f^*c_1(\FCC).$$
\end{prop}
\dem La preuve de notre proposition d{\'e}coule directement du diagramme
commutatif suivant :
$$\xymatrix{(X/S)_{cris}^m \ar[r]^{f_{cris}^m}
  \ar[d]_{{u_{X/S}}} & (Y/S)_{cris}^m \ar[d]^{u_{Y/S}} \\
X_{Zar} \ar[r]_{f^m} & Y_{Zar}}$$
et du fait que les suites exactes utilis\'ees sont fonctorielles.

\findem

\label{m-def_cc}

Maintenant que l'on a construit la premi{\`e}re classe de Chern d'un
faisceau inversible, nous allons en d{\'e}duire les autres par la m{\'e}thode
habituelle \cite{Gro3}.

On se donne $\ECC$ un faisceau localement libre de rang $r$ sur $X$,
on note $\PM=\PM(\ECC)$ son fibr{\'e} projectif et $\OCC_{\PM}(1)$ le
faisceau canonique sur ce dernier. Le morphisme structural de $\PM$
sera not{\'e} $p : \PM \to X$. 

On a alors le diagramme commutatif de morphismes de topos suivant :
$$\xymatrix{(\PM/S)_{cris}^m \ar[r]^{p_{cris}^m}
  \ar[d]_{{u_{\PM/S}}} & (X/S)_{cris}^m \ar[d]^{u_{X/S}} \\
\PM_{Zar} \ar[r]_{p} & X_{Zar}.}$$

On note $$\xi = c_{1,m}(\OCC_{\PM}(1)) \in H^2_{m-cris}(\PM).$$
On peut alors voir pour tout $i$, $\xi^i$ comme un morphisme dans
$D^+((X/S)^m_{cris},Ab)$ :
$$\xi^i : \ZM_X \to \RM p^m_{cris*}\OCC_{\PM/S}^m[2i].$$

De m{\^e}me, le morphisme $p^m_{cris}$ nous donne un morphisme :
$$\OCC_{X/S}^m \to \RM p_{cris*}^m \OCC_{\PM/S}^m.$$

On en d{\'e}duit alors comme dans le cas classique :
\begin{equation}
\label{m-oplus}
\bigoplus_{i=0}^{r-1} \xi^i : \bigoplus_{i=0}^{r-1} \OCC_{X/S}^m[-2i] \to \RM p^m_{cris*}\OCC_{\PM/S}^m.
\end{equation}

On a alors la proposition suivante :

\begin{prop}
\label{m-structure}
Avec les notations pr{\'e}c{\'e}dentes, l'application (\ref{m-oplus}) est un isomorphisme.
\end{prop}

\dem
On remarque que, la question {\'e}tant locale, on peut supposer que
$\ECC=\OCC_{X}^r$ et qu'il existe une immersion ferm{\'e}e de
$X$ dans $Y$ un $S$-sch{\'e}ma formel lisse. En notant $\widehat{\PG}$
l'enveloppe {\`a} puissance $m$-divis{\'e}e de cette immersion on a le
diagramme :
$$\xymatrix{\PM^r_{X} \ar[r] \ar[d] & \PM^r_{\widehat{\PG}} \ar[d]
  \\
            X \ar[r]   & \widehat{\PG} }$$
Dans cette situation, la cohomologie $m$-cristalline se calcule comme la cohomologie de De Rham de l'enveloppe \`a puissances divis\'ees. Le faisceau $\OCC_{\PM^r_{X}}(1)$ se relevant en
  $\OCC_{\PM^r_{\widehat{\PG}}}(1)$, on est ramen\'e au m\^eme probl\`eme pour la cohomologie de De Rham de $\PM^r$ sur $\widehat{\PG}$. La proposition est connue dans ce cas.

\findem

\begin{coro}
\label{m-coroll}
Avec les notations ci-dessus on a pour tout $n$ la
d{\'e}composition suivante :
$$H^n_{m-cris}(\PM)=\bigoplus_{i=0}^{r-1}H^{n-2i}_{m-cris}(X).$$

\end{coro}

On d{\'e}finit les autres classes de Chern comme dans le cas
classique. On applique la d{\'e}composition du corollaire \ref{m-coroll} {\`a}
$\xi^r$ et on obtient :
$$\xi^r=\sum_{i=1}^r (-1)^{i+1} c_{i,m}(\ECC)\xi^{r-i},$$
avec :
$$c_{i,m}(\ECC) \in H^{2i}_{m-cris}(X).$$
Par d{\'e}finition $c_i(\ECC)$ est la $i$-{\`e}me classe de Chern de $\ECC$.

\rqe Comme dans la section \ref{propre}, on note $c_0(\ECC) =1$ dans $H^0_{m-cris}(X)$.

Nous allons nous int\'eresser \`a l'additivit\'e des classes de Chern d\'efinies ci-dessus. La d\'emonstration utilisant une comparaison avec la cohomologie de De Rham, nous nous contenterons de regarder la situation apr\`es tensorisation par $\QM$ \footnote{Soient $X$ une $k$-vari\'et\'e et $X^{(m)}$ le pull-back par le $m$-ieme it\'er\'e du Frobenius de $k$. Pour tout faisceau localement libre $\ECC$ sur $X$, on notera $\ECC^m$ sont pull-back sur $X^{(m)}$. Il semble clair alors que la classe de Chern $c_{i,m}(\ECC)$ \`a valeurs dans $H^{2i}_{m}(X)$ soit l'image de la classe $c_{i,0}(\ECC^{(m)})$ par le morphisme 
$$F_m^* : H^{2i}_{cris}(X^{(m)}) \to  H^{2i}_{m-cris}(X^),$$
dont on trouve la d\'efinition dans \cite{Berth6}. On d\'eduit alors de l'additivit\'e des classes de Chern cristallines de niveau $0$, l'additivit\'e des classes de Chern $m$-cristallines.}.

Plus pr\'ecis\'ement, on note 
$$H^i_{m,\QM}(X) := H^i((X/S)^m_{cris}, \OCC_{X/S}\otimes \QM)$$
et, pour tout faisceau localement libre $\ECC$ sur $X$
$$c_{i,m}^{\QM}(\ECC) \in H^{2i}_{m,\QM}(X),$$
la $i$-eme classe de Chern d\'efinie comme ci-dessus.

\begin{prop}[Additivit\'ee]
Soient $X$ une $k$-vari\'et\'e, $\ECC$, $\ECC'$ et $\ECC''$ trois faisceaux localement libres tels que l'on ait la suite exacte :
$$0 \to \ECC' \to \ECC \to \ECC'' \to 0.$$
On a alors pour tout $m$ et tout $i$ on a dans $H^{2i}_{m,\QM}(X)$ :
$$c_{i,m}^\QM(\ECC) = \sum_{j=0}^ic_{j,m}^\QM(\ECC')c_{i-j,m}^\QM(\ECC'').$$
\end{prop}

\dem
 La fonctorialit\'e de nos constructions permet de les g\'en\'eraliser au cas des vari\'et\'es simpliciales. On peut donc appliquer la m\'ethode classique \cite[pp 217-224]{Gill} qui consiste \`a se ramener \`a une situation universelle. Remarquons que notre th\'eorie cohomologique ne v\'erifie pas les axiomes demand\'es par {\it loc. cit.}, cependant poour ramener l'additivit\'e des classes de Chern \`a la situation universelle sur $B_{\bullet}GL$ on n'a besoin que du corollaire \ref{m-coroll}. A partir de l\`a, on ne regarde que des vari\'et\'es lisses relevables. On d\'eduit donc notre proposition du r\'esultat similaire en cohomologie de De Rham.

\findem

\begin{prop}
Soient $X$ une $k$-vari\'et\'e et $\ECC$ un faisceau localement libre sur $X$ . Pour tous $m, m'$ avec $m \leq m'$, on a pour tout $i$ :
$$\phi_{m',m} c_{i,m'}(\ECC) = p^{i(m'-m)
}c_{i,m}(\ECC)$$
o\`u $\phi_{m',m} : H^{2i}((X/S)_{cris}^{m'},\OCC_{X/S}^{m'}) \to  H^{2i}((X/S)_{cris}^{m},\OCC_{X/S}^{m})$ est le morphisme induit par $u_{m',m}$.
\end{prop}

\dem Il suffit de le voir pour la premi\`ere classe de Chern d'un faisceau inversible. Ce dernier cas d\'ecoule directement de la d\'efinition.

\findem

\subsection{Le topos $(X/S)_{cris}^\bullet$}

Nous allons construire un topos permettant de calculer la cohomologie rigide de mani{\`e}re cristalline. Ce topos se construit comme le topos associ\'e \`a un diagramme de topos \cite[IV,7]{SGA4} ou \cite[VI]{Ill2}. Le cas particulier nous int\'eressant, \`a savoir le cas des syst\`emes inductifs de topos est repris dans \cite[II]{moi}.

\begin{Defi}
Avec les notations pr{\'e}c{\'e}dentes, on notera $(X/S)_{cris}^\bullet$ le
topos associ{\'e} au diagramme de topos
$((X/S)_{cris}^m, i_{m',m})_{m \in \NM}$.
\end{Defi}

\notation On notera souvent un objet de ce topos de la mani{\`e}re
suivante
$$\FCC^\bullet = (\FCC^0 \leftarrow \cdots \leftarrow \FCC^m
\leftarrow \cdots).$$
On pose alors :
$$\OCC_{X/S}^\bullet = (\OCC_{X/S}^0 \leftarrow \cdots \leftarrow \OCC_{X/S}^m
\leftarrow \cdots),$$
les fl{\`e}ches de transition {\'e}tant, pour $m \leq m'$,  les morphismes $u_{m',m}$
d{\'e}finis pr{\'e}c{\'e}demment.
De la m{\^e}me mani{\`e}re, on pose :
$$\ICC_{X/S}^\bullet = (\ICC_{X/S}^0 \leftarrow \cdots \leftarrow \ICC_{X/S}^m
\leftarrow \cdots)$$

On consid{\`e}re alors le topos $(X/S)_{cris}^\bullet$ comme annel{\'e}
par $\OCC_{X/S}^\bullet$. On a les morphismes de topos annel\'es suivants :
$$p_m : ((X/S)_{cris}^m, \OCC_{X/S}^m) \to ((X/S)_{cris}^\bullet, \OCC_{X/S}^\bullet).$$

Le lemme \ref{commut} nous dit que les morphismes de topos 
$$u_{X/S}^m : (X/S)_{cris}^m \to X_{Zar} \hbox{ (resp. }  i_{X/S}^m : X_{Zar} \to (X/S)_{cris}^m)$$
permettent de construire un morphisme de topos :
$${u^\bullet_{X/S}} : (X/S)_{cris}^\bullet \to
X_{Zar}^\bullet \hbox{ (resp. } {i^\bullet_{X/S}} : X_{Zar}^\bullet \to (X/S)_{cris}^\bullet).$$

De m{\^e}me, les foncteurs sections globales 
$\Gamma((X/S)^m_{cris},-) : (X/S)^m_{cris} \to (Ens)$ donnent naissance {\`a} un foncteur :
$$\Gamma^\bullet((X/S)_{cris}^\bullet, -) :
(X/S)_{cris}^\bullet \to (Ens)^\bullet.$$

De plus, l'{\'e}galit{\'e} de foncteurs :
$$ \Gamma((X/S)_{cris}^m, -)= \Gamma (X_{Zar},
-) \circ 
u^m_{X/S*} $$
implique que :
$$\Gamma^\bullet((X/S)_{cris}^\bullet, -)= \Gamma^\bullet (X_{Zar}^\bullet,
-) \circ 
u^\bullet_{X/S*} $$
o{\`u} $\Gamma^\bullet (X_{Zar}^\bullet,
-)$ est le foncteur de
$X_{Zar}^\bullet$ dans $(Ens)^\bullet$ d{\'e}fini de la
m{\^e}me mani{\`e}re que $\Gamma^\bullet((X/S)_{cris}^\bullet, -)$.

\vspace{.3cm}

\rqe On fera attention {\`a} ne pas confondre les foncteurs $\Gamma^\bullet((X/S)_{cris}^\bullet, -)$ et $\Gamma^\bullet (X_{Zar}^\bullet,
-)$ avec les foncteurs
sections globales des topos $(X/S)_{cris}^\bullet$ et $X_{Zar}^\bullet$. On a cependant
$$\Gamma((X/S)_{cris}^\bullet, - ) = \limpro \circ \Gamma^\bullet((X/S)_{cris}^\bullet, -) \hbox{ et } \Gamma(X_{Zar}^\bullet, - ) = \limpro \circ \Gamma^\bullet(X_{Zar}^\bullet, -).$$
\notation Par la suite, nous noterons $\RM \u{\Gamma}_{cris\bullet}(X)$ le complexe $\Rlimproj{m}((\RM {u_{X/S}^\bullet}_* \OCC_{X/S}^\bullet)\otimes \QM)$ et 
$$H^*_{cris\bullet}(X) := H^*(X_{Zar}^\bullet, (\RM {u_{X/S}^\bullet}_* \OCC_{X/S}^\bullet)\otimes \QM)= H^*(X_{Zar},\RM \u{\Gamma}_{cris\bullet}(X)) .$$

\subsection{Interpr{\'e}tation cristalline de la cohomologie rigide}
Les d{\'e}finitions pr{\'e}c{\'e}dentes permettent d'\'enoncer l'interpr{\'e}tation cristalline de la cohomologie rigide.

\begin{theo}[Berthelot]
\label{berthelot}
Soit $X$ une $k$-vari{\'e}t{\'e} propre. En reprenant les notations pr\'ec\'edentes on a
$$\RM \u{\Gamma}_{rig}(X/K) \iso \RM \u{\Gamma}_{cris\bullet}(X).$$
\end{theo}
\dem Remarquons juste que la d\'emonstration de ce th\'eor\`eme utilise un calcul de la cohomologie $m$-cristalline \`a l'aide de la cohomologie de De Rham qui ne se fait qu'a torsion pr\'es, il est donc n\'ecessaire de tensoriser par $\QM$. La d\'emonstration se trouve dans \cite{Berth6}.
\findem

En passant \`a la cohomologie on obtient :
$$H^i_{rig}(X/K) \iso H^i_{cris\bullet}(X).$$

\subsection{Construction des classes de Chern sur $(X/S)^\bullet_{cris}$}
La construction est similaire \`a celle des classes de Chern $m$-cristalline.

Avec les notations pr\'ec\'edentes, on a la
suite exacte dans $(X/S)_{cris}^\bullet$ : 
$$0 \to 1 + \ICC_{X/S}^{\bullet} \to (\OCC_{X/S}^*)^\bullet \to
i_*^\bullet(\OCC_{X}^*) \to 0.$$
Pour le d{\'e}montrer il suffit d'utiliser la famille conservative des $p^{-1}_m$.
On obtient donc un
morphisme dans $D^+((X/S)^\bullet_{cris},Ab)$ :
$$i^\bullet _*(\OCC_{X}^*)[1] \to (1+\ICC_{X/S}^{\bullet})[2].$$
En projetant ce dernier sur le topos zariskien on obtient dans
$D^+(X_{Zar}^\bullet,Ab)$ le morphisme :

\begin{equation}
\label{prem}
(\OCC_{X}^*)^\bullet[1] \to \RM u^\bullet_{X/S*}(1+\ICC^\bullet_{X/S})[2].
\end{equation}

Cependant, les morphismes $\psi_m$ d\'efinis pr\'ec\'edemment ne construisent pas un morphisme sur
les syst{\`e}mes projectifs car ils ne commutent pas aux fl{\`e}ches des
syst{\`e}mes projectifs.

On consid{\`e}re le syst{\`e}me projectif :
$$(1+\ICC^\bullet_{X/S})(1)=(1+\ICC^0_{X/S} \gsurfleche{p}
\cdots \gsurfleche{p} 1+\ICC^m_{X/S} \gsurfleche{p} \cdots )$$
o{\`u} les fl{\`e}ches de transition sont l'{\'e}l{\'e}vation {\`a} la puissance $p$ des
fl{\`e}ches classiques.
 
Les morphismes $\psi_m$ d{\'e}finissent alors :
$$\psi^\bullet : (1+\ICC^\bullet_{X/S})(1) \to
\JCC^\bullet_{X/S}.$$

Apr{\`e}s projection sur le topos zariskien, on obtient :

$$\RM
u^\bullet_{X/S*}((1+\ICC^\bullet_{X/S})(1))[2] \to \RM
u^\bullet_{X/S*}(\JCC^\bullet_{X/S})[2].$$
En utilisant composant cette derni\`ere avec la projection du  morphisme $\JCC^\bullet_{X/S} \to \OCC_{X/S}^\bullet$, on a :

\begin{equation}
\label{deux}
\RM
u^\bullet_{X/S*}((1+\ICC^\bullet_{X/S})(1))[2] \to \RM
u^\bullet_{X/S*}(\OCC_{X/S}^\bullet)[2].
\end{equation}
 
On peut d'autre part construire un  morphisme de 
$(1+\ICC^\bullet_{X/S})(1)$ dans $1+\ICC^\bullet_{X/S}$ d{\'e}fini
par l'{\'e}l{\'e}vation {\`a} la puissance $p^{m}$ en {\it degr{\'e}}
$m$.

On remarque alors que ce morphisme induit apr{\`e}s tensorisation par
$\QM$ un isomorphisme dans $D^+(X_{Zar}^\bullet, Ab)$ :
\begin{equation}
\label{trois}
\RM u^\bullet_{X/S*}((1+\ICC^\bullet_{X/S})(1))\otimes \QM \iso \RM u^\bullet_{X/S*}(1+\ICC^\bullet_{X/S})\otimes \QM.
\end{equation}

En composant \ref{prem}, l'inverse de \ref{trois} et \ref{deux}, on obtient, apr\`es passage \`a la limite :

\begin{equation}
\label{c1}
c_{1,\bullet} : (\OCC_{X}^*)[1] \to \RM \u{\Gamma}_{cris\bullet}(X)[2].
\end{equation}

Par suite en passant \`a la cohomologie on a :
\begin{equation}
\label{c1g}
c_{1,\bullet} : H^1(X,\OCC_X^*) \to H^2_{cris\bullet}(X).
\end{equation}

Comme dans le cas $m$-cristallin, la premi\`ere classe de Chern permet de calculer la cohomologie des fibr\'es projectifs. On peut donc construire classiquement les autres classes de Chern.

\rqe Comme dans la section \ref{propre}, on note $c_0(\ECC) =1$ dans $H^0_{cris\bullet}(X)$. De plus, si $\ECC$ est un faisceau inversible, on retrouve la premi\`ere classe de Chern construite auparavant.

\begin{prop}
Les classes de Chern $c_{i,\bullet}$ v\'erifient les propri\'et\'es d'additivit\'ee et de fonctorialit\'e.
\end{prop}

\dem Cela se d\'emontre comme les \'enonc\'es similaires pour la cohomologie $m$-cristalline.

\findem

\subsection{Comparaison avec les classes de Chern rigides}
Nous allons voir que les classes de Chern que nous venons de
construire co\"incident avec les classes de Chern rigides d{\'e}finies \`a la section  \ref{propre}. On se donne $X$, une $k$-vari\'et\'e propre, ainsi qu'un plongement $X \hookrightarrow \YCC$ dans un sch\'ema formel sur $S=\spf{C}$ lisse au voisinage de $X$.
\begin{prop}
On a le diagramme commutatif suivant :
$$\xymatrix{ & H^2_{cris\bullet}(X)\\
H^1(X, \OCC_X^*) \ar[ur]^{c_{1,\bullet}} \ar[dr]_{c_{1,rig}} & \\
& H^2_{rig}(X/K). \ar[uu]_{\wr} }$$
o{\`u} $c_{1,\bullet}$ est le morphisme \ref{c1g} et la fl\`eche verticale se d\'eduit de celle de \ref{berthelot}.
\end{prop}
\dem

\begin{coro}
Les classes de Chern rigides v\'erifient la propri\'et\'es d'additivit\'e.
\end{coro}

\label{112}

Nous allons, comme dans \cite{Berth-Ogu2}, r{\'e}interpr\'eter la construction de nos classes de Chern cristallines {\`a} l'aide de faisceaux Zariskiens. Nous utiliserons alors des r{\'e}solutions de Cech afin de les comparer avec la d{\'e}finition de la section \ref{propre}. 

Dans $D^+(X_{Zar})$ nous noterons avec un $\QM$ en indice les
complexes obtenus apr{\`e}s tensorisation par $\QM$.

On notera $Y_n$ la r\'eduction modulo $p^{n+1}$ de $\YCC$, pour tout $m$ et tout $n$, on note alors $\PG^m_n$ l'enveloppe {\`a} puissance $m$-divis{\'e}e de l'immersion de $X$ dans $Y_n$, $\hat{\PG}^m$ sa compl{\'e}tion $p$-adique et $\hat{\PG}^\bullet$ le syst{\`e}me projectif associ{\'e}.

On
regarde le complexe de De Rham et le complexe de De Rham multiplicatif
dans $X_{Zar}^\bullet$, {\`a} savoir :
$$\Omega^\star_{\widehat{\PG}^\bullet/S} = \OCC_{\widehat{\PG}^\bullet} \surfleche{d} \Omega^1_{\widehat{\PG}^\bullet/S} \to
\Omega^2_{\widehat{\PG}^\bullet/S} \to \cdots $$
et :
$$\Omega^\times_{\widehat{\PG}^\bullet/S} = \OCC_{\widehat{\PG}^\bullet}^* \surfleche{dlog} \Omega^1_{\widehat{\PG}^\bullet/S} \to
\Omega^2_{\widehat{\PG}^\bullet/S} \to \cdots $$
et on note $\KCC_{\widehat{\PG}^\bullet/S}^\times$ (resp. $\ICC_{\widehat{\PG}^\bullet/S}^\star$) le
noyau de l'application de complexes $\Omega^\times_{\widehat{\PG}^\bullet/S} \to
(\OCC_{X}^*)^\bullet$
(resp. $\Omega^\star_{\widehat{\PG}^\bullet/S} \to
(\OCC_{X})^\bullet$).
On rappelle que $\ICC^{m}$ est l'id{\'e}al de $X$ dans $\widehat{\PG}^{m}$.

De plus, comme pr{\'e}c{\'e}demment, on note :
$$(1+\ICC^\bullet)(1) = (1+\ICC^0 \gsurfleche{p} 1+\ICC^1
\gsurfleche{p} \cdots \gsurfleche{p} 1+\ICC^m \gsurfleche{p} \cdots)$$
et :
$$\Omega_{\widehat{\PG}^\bullet/S}^i(1) =
(\Omega_{\widehat{\PG}^0/S}^i\gsurfleche{p}\Omega_{\widehat{\PG}^1/S}^i\gsurfleche{p} \cdots \gsurfleche{p}  \Omega_{\widehat{\PG}^m/S}^i\gsurfleche{p} \cdots).$$
Enfin on pose :
$$\KCC^\times_{\widehat{\PG}^\bullet/S}(1) = (1+\ICC^\bullet)(1)
\surfleche{dlog} \Omega_{\widehat{\PG}^\bullet/S}^1(1) \to
\Omega_{\widehat{\PG}^\bullet/S}^2(1) \to \cdots $$

On a alors un morphisme :
$$(\KCC^\times_{\widehat{\PG}^\bullet/S}(1))_{\QM} \to (\Omega^\star_{\widehat{\PG}^\bullet/S})_{\QM}$$
d\'efini par les applications :
$$p^m\log : (1+\ICC^m) \to \OCC_{\widehat{\PG}^m}$$
et 
$$\times p^m : \Omega^i_{\widehat{\PG}^m/S} \to \Omega^i_{\widehat{\PG}^m/S}.$$

Avec ces notations on a le diagramme commutatif suivant dans $D^+(X^\bullet_{zar})$ :

$$\xymatrix{\RM u^{\bullet}_{X/S*}(i^\bullet_*\OCC^*_{X})[1]
  \ar[r]^{\partial} \ar@{=}[d]
& \RM u^{\bullet}_{X/S*}(1 +
  \ICC^{\bullet}_{X/S})_{\QM}[2] \ar[d]_{\wr} 
&\RM u^{\bullet}_{X/S*}((1 +
  \ICC^{\bullet}_{X/S})(1))_{\QM}[2] \ar[l]_\sim \ar[d]_{\wr} \ar[r]&
\RM u^{\bullet}_{X/S*}(\OCC^{\bullet}_{X/S})_\QM[2] \ar[d]_{\wr}
\\
(\OCC_{X}^*)^\bullet[1] \ar[r]^{\partial} 
& (\KCC^\times_{\widehat{\PG}^\bullet/S})_{\QM}[2] 
  &(\KCC^\times_{\widehat{\PG}^\bullet/S}(1))_{\QM}[2] \ar[l]_\sim \ar[r]
& (\Omega^\star_{\widehat{\PG}^\bullet/S})_{\QM}[2]}$$

La premi{\`e}re classe de Chern \'etant, par d\'efinition, l'application induite sur la cohomologie par la premi\`ere ligne du diagramme, elle peut \^etre calcul{\'e}e par la deuxi{\`e}me. Nous allons maintenant en faire un calcul explicite {\`a} l'aide de cocycles de Cech.


On sait que pour tout $\sigma \in H^1(X,\OCC_X^*)$, il existe un
recouvrement $\UG$ de $\YCC$ tel que $\sigma$ soit l'image d'un {\'e}l{\'e}ment
(encore not{\'e} $\sigma$) de $H^1(\UG_X, \OCC_X^*)$. 

\vspace{.3cm}

\notation on reprend les notations de la section \ref{propre}.
On notera de plus, $\ICC^m_{i}$ l'id{\'e}al de $X_{i}$ dans $\widehat{\PG}^m_{i}$ et
$I^m_{i}=\Gamma(\widehat{\PG}^m_{i},\ICC^m_{i})$. De
m{\^e}me, on note $\widehat{\PG}^m_i = \spf{\widehat{P}^m_{i}}$.  En notant $P^m_{n,i}$ l'enveloppe \`a puissances divis\'ees de niveau $m$ de $(\ACC_i/p^{n+1}\ACC_i, I_i/p^{n+1}I_i)$ o\`u $I_i$ est l'id\'eal de $A_i$ dans $\ACC_i$, on a 
$$\widehat{P}^m_{i} = \limproj{n} P^m_{n,i}.$$
On a alors des morphismes : 
$$\pi^m_{i} : \widehat{P}^m_{i} \to A_i.$$
On notera pour tout $m'>m$, $\varepsilon^{m,m'}_{i}$ la fl\`eche entre $\widehat{P}^{m'}_{i}$ et $\widehat{P}^{m}_{i}$. Dans le cas particulier o\`u $m' = m+1$ nous la noterons simplement $\varepsilon^m_{i}$.

Avec ces notations, $\sigma$ est donc repr{\'e}sent{\'e} par un cocycle $u_{ij} \in C^1(\UG_X, \OCC_X^*)$. On regarde donc le cocycle $u^m_{ij} \in C^1(\UG, (\OCC_X^*)^\bullet)$ d{\'e}fini par $u^m_{ij}=u_{ij}$ pour tout $m$. Ce cocycle d\'efinit une application dans $D^+(X_{Zar})$ : $\ZM \to (\OCC_{X}^*)^\bullet[1].$
Nous allons montrer que l'application de $D^+(X_{Zar})$ 
$$\ZM \to (\Omega^\star_{\widehat{\PG}^\bullet/S})_{\QM}[2]$$
obtenue comme composition de cette application avec la deuxi\`eme ligne du diagramme pr\'ec\'edent est repr\'esentable par un cocycle de $C^2(\UG, \Omega^{\star}_{\widehat{\PG}^\bullet/S}\otimes  \QM).$

\begin{Defi}
Avec les notations pr{\'e}c{\'e}dentes, pour tout $x \in A_i$ on appelle
rel{\`e}vement compatible de $x$ une famille $(\widetilde{x}^m) \in \Pi_{m\in \NM}
\widehat{P}^{m}_i$ telle que :
\begin{itemize}
\item pour tout $m$ on ait $\pi^m_i (\widetilde{x}^m) = x$
\item pour tout $m'>m$ on ait $\varepsilon^{m,m'}_i(\widetilde{x}^{m'})
  = \widetilde{x}^m.$
\end{itemize}
\end{Defi}
Le diagramme commutatif suivant :
$$\xymatrix{  & \widehat{P}^{m'}_i \ar[d]^{\varepsilon^{m,m'}_i}
  \ar@/^1.5pc/[ddr]^{\pi^{m'}_i} & \\
    & \widehat{P}^{m}_i \ar[dr]^{\pi^{m}_i} & \\
\ACC_i \ar[rr] \ar@/^1.5pc/[uur] \ar[ur] & & A_i}$$
montre que pour construire un rel\`evement compatible, il suffit de relever $x$ dans $\ACC_i$ et d'utiliser les morphismes $\ACC_i \to \widehat{P}^{m}_i.$ En particulier, il existe toujours des rel\`evements compatibles.

On appellera syst{\`e}me compatible un {\'e}l{\'e}ment v{\'e}rifiant la
deuxi{\`e}me condition.

Soit $(\widetilde{u}^m_{ij}) \in \Pi_{m\in \NM}\widehat{P}^{m}_{ij}$ un rel{\`e}vement compatible de $u_{ij}=u^m_{ij}$.
\begin{lemme}
Avec les notations pr{\'e}c{\'e}dentes, $\widetilde{u}^{m}_{ij}$ est
inversible dans $\widehat{P}^{m}_{ij}$.
\end{lemme}
\dem
On rel{\`e}ve aussi $(u_{ij})^{-1}$ de mani{\`e}re compatible en $\widetilde{v}^m_{ij}$. On a alors :
$$\widetilde{u}^{m}_{ij}.\widetilde{v}^{m}_{ij}=1+x^{m}_{ij} \in 1+I^{m}_{ij}.$$
De plus, $I^m_{ij}$ {\'e}tant un $m$-PD-id{\'e}al, pour tout $x \in I^{m}_{ij}$
et tout $k\in \NM$, on a :
$$x^k = q_k! x^{\{k\}_{m}}$$ 
o{\`u} $k=p^mq_k +r_k$ et $0 \leq r_k < p^m$.
A partir de l{\`a}, on voit que $x^k$ tend $p$-adiquement vers $0$ quand
$k$ tend vers l'infini. Comme $\widehat{P}^m_{ij}$ est complet,
la s\'erie :
$$\sum_{k=0}^\infty (-x^m_{ij})^k$$
converge vers $(1+x^m_{ij})^{-1}$. L'inverse de $\widetilde{u}^{m}_{ij}$
est donc :
$$\widetilde{v}^{m}_{ij}.\sum_{k=0}^\infty(-x^m_{ij})^k.$$

\findem

\rqe ce lemme est l'anologue $m$-cristallin du lemme \ref{laumon}.

De plus il est clair d'apr{\`e}s la d{\'e}finition que
$\widetilde{x}^{m}_{ij}$ et donc $(\widetilde{u}^{m}_{ij})^{-1}$ sont des
syst{\`e}mes compatibles.
 A partir de l{\`a}, on consid{\`e}re :
$$\dlog{\widetilde{u}^{m}_{ij}} \in C^1(\UG, \Omega^1_{\widehat{\PG}^m/S})$$
et :

$$h^m_{ijk}=\widetilde{u}^{m}_{ij}.(\widetilde{u}^{m}_{ik})^{-1}.\widetilde{u}^{m}_{jk}\in
C^2(\UG,1+\ICC^m).$$
Ils d{\'e}finissent donc un {\'e}l{\'e}ment de $C^2(\UG, \KCC^\times_{\widehat{\PG}^\bullet/S})$. Des calculs similaires \`a ceux de la section \ref{propre} montrent que si on modifie les rel\`evements choisis ou le choix du cocycle, on modifie cette classe d'un cobord. On obtient donc bien une application dans la cat\'egorie d\'eriv\'ee :
$$\ZM \to \KCC^\times_{\widehat{\PG}^\bullet/S}[2].$$

On inverse alors, apr{\`e}s tensorisation par $\QM$, le morphisme $\KCC^\times_{\widehat{\PG}^\bullet/S}\otimes \QM \to \KCC^\times_{\widehat{\PG}^\bullet/S}(1)\otimes \QM$ en consid{\'e}rant le morphisme qui est d\'efini en degr{\'e} $m$ par la division par $p^m$. On obtient donc le cocycle :
$$\left(\dlog{\widetilde{u}^m_{ij}}\otimes \frac{1}{p^m}, \widetilde{u}^m_{ij}.(\widetilde{u}^m_{ik})^{-1}.\widetilde{u}^m_{jk}\otimes \frac{1}{p^m}\right) \in  C^2(\UG, \KCC^\times_{\widehat{\PG}^\bullet/S}(1)\otimes \QM).$$
On conclut alors notre calcul en utilisant l'application {\'e}gale en degr\'e $m$ \`a  la multiplication par $p^m$ sur la premi\`ere composante et $p^m\log$ sur la deuxi\`eme. Ce qui nous donne :
$$\left(\dlog{\widetilde{u}^m_{ij}}\otimes 1, \log(\widetilde{u}^m_{ij}.(\widetilde{u}^m_{ik})^{-1}.\widetilde{u}^m_{jk})^{p^m}\otimes \frac{1}{p^m}\right) \in  C^2(\UG, \Omega^\star_{\widehat{\PG}^\bullet/S})_{\QM}.$$

Ce cocycle correspond, par le morphisme du th{\'e}or{\`e}me \ref{berthelot}, {\`a} celui qui a {\'e}t{\'e} construit \`a la section\ref{propre}.
\findem

\begin{coro}
Soit $X$ une $k$-vari\'et\'e propre, on a le diagramme commutatif suivant :
$$\xymatrix{ & H^2((X/S)_{cris}, \OCC_{X/S}) \otimes \QM \\ 
H^1(X, \OCC_X^*) \ar[ur]^{c_{1,cris}} \ar[dr]_{c_{1,rig}} &  \\
  & H^2_{rig}(X/K) \ar[uu] }$$
o\`u la fl\`eche $c_{1,cris}$ est la premi\`ere classe de Chern cristalline de \cite{Berth-Ill} et la fl\`eche verticale est celle qui est donn\'ee par
$$H^2_{rig}(X/K) \to H^2_{cris\bullet}(X) \to  H^2((X/S)_{cris}, \OCC_{X/S}) \otimes \QM.$$
\end{coro}

\dem
Gr\^ace \`a la proposition, il suffit de comparer $c_1^\bullet$ avec $c_{1,cris}$. Or cette derni\`ere s'obtient \`a partir de la premi\`ere en appliquant le foncteur $p_0^*$.
\findem

\begin{coro}
\label{conv}
On a 
$$c_{i,rig}(\ECC) = c_{i,conv}(\ECC)$$
via l'isomorphisme de \cite[0.6.7]{Ogu}
$$H^{2i}_{rig}(X/K) \iso H^{2i}((X/S)_{conv}, \KCC_{X/S})$$
o\`u $c_{i,conv}$ sont les classes de Chern d\'efinies dans \cite[A]{Niz}.
\end{coro}

\dem En notant $T = (X,\YCC,Id_X)$ on a les isomorphismes :
$$\Rlimproj{m}(\RM u^\bullet_{X/S*}(\OCC_{X/S}^\bullet) \otimes {\QM}) \iso \RM sp_*\Omega^\star_{]X[} \iso (\KCC_{X/S})_T.$$

D\`es lors la construction du $c_{1,conv}$ de \cite[A]{Niz} est la m\^eme que celle du $c_{1,\bullet}$ ci dessus. Les d\'ecompositions de $H^*_{rig}(\PM_X/K)$ comme  $H^*_{rig}(X/K)$-module et $H^*((\PM_X/S)_{conv}, \KCC_{\PM_X/S})$ comme $H^*((X/S)_{conv}, \KCC_{X/S})$-module \'etant compatibles, notre comparaison se prolonge aux $c_i$. 

\findem

\section{Cas des vari\'et\'es ouvertes}
Nous allons nous int\'eresser au cas g\'en\'eral.

\subsection{M{\'e}thode g{\'e}n{\'e}rale - Les th{\'e}or{\`e}mes de prolongements}
Nous allons construire la premi{\`e}re classe de Chern d'un faisceau
inversible en montrant qu'il existe une compactification de $X$ sur
laquelle on peut prolonger ce faisceau en un faisceau inversible. On pourra alors
conclure en montrant que la classe obtenue, en revenant, par fonctorialit\'e, dans la
cohomologie rigide de $X$, ne d{\'e}pend pas des choix faits.

L'ingr\'edient principal pour montrer l'existence de ces prolongements de faisceaux est le th\'eor\`eme de platification par \'eclatement de Raynaud-Gruson \cite{Gru-Ray}. Nous allons avoir besoin de terminologie. Soient $X$ une $k$-vari\'et\'e, $Y$ une sous-vari\'et\'e ferm\'ee et $\FCC$ un $\OCC_X$-module. Si on note $\pi : X' \to X$ l'\'eclat\'e de $Y$ dans $X$, on appellera transform\'e strict de  $\FCC$, et on notera TS$(\FCC)$, le quotient de $\pi^*(\FCC)$ par le sous-faisceau engendr\'e par les sections \`a support dans $\pi^{-1}(Y).$
Il est direct de montrer que le transform\'e strict d'un faisceau coh\'erent (resp. localement libre) est coh\'erent (resp. localement libre).

Passons \`a la d\'emonstration des lemmes.

\begin{lemme} 
\label{lemme1}
Soit $X$ une $k$-vari{\'e}t{\'e}, il existe une
  $k$-vari{\'e}t{\'e} propre $\ov{X}$ et une immersion ouverte $j : X
  \hookrightarrow \ov{X}$ telles que $Z=\ov{X}-X$ soit un diviseur.
\end{lemme}
\dem D'apr{\`e}s Nagata \cite{Nag}, il existe une $k$-vari{\'e}t{\'e} propre
$\ov{X}_1$ et une immersion ouverte $j_1 : X \hookrightarrow
\ov{X}_1$. Si on note $Z_1 = \ov{X}_1 - X$, on consid{\`e}re
l'\'eclatement $\ov{X}$ de $\ov{X}_1$ le long de $Z_1$. C'est encore une
$k$-vari{\'e}t{\'e} propre. De plus, l'\'eclatement ne modifiant pas $\ov{X}_1 - Z_1$, $j_1$ se prolonge en une immersion
ouverte $j : X \hookrightarrow \ov{X}$ et $Z = \ov{X} - X$ est un diviseur.
\findem

\rqe Par la suite nous noterons $(\ov{X},j)$ la donn\'ee consistant en la vari\'et\'e propre $\ov{X}$ et l'immersion ouverte $j : X \hookrightarrow \ov{X}$. Nous supposerons toujours, sauf mention explicite que nos compactifications v\'erifient la condition du lemme \ref{lemme1}.

\begin{theo}
\label{theo1}
Soient $X$ une $k$-vari{\'e}t{\'e} et $\ECC$ un faisceau localement libre
de rang $r$ sur $X$. Il existe une compactification de $X$, $j : X
\hookrightarrow \ov{X}$
et $\ov{\ECC}$ un faisceau localement libre de rang $r$ sur $\ov{X}$
tel que :
$$j^*(\ov{\ECC})=\ECC.$$
\end{theo}
\dem On se donne une compactification $(\ov{X}_1, j_1)$ de
$X$. D'apr{\`e}s \cite[6.9.8]{EGA1N}, on sait qu'il existe $\ECC_1$ un
sous-faisceau coh{\'e}rent de ${j_1}_*(\ECC)$ tel que
$j_1^*(\ECC_1)=\ECC$.  On sait alors par le th{\'e}or{\`e}me de
platification par {\'e}clatement de Raynaud-Gruson \cite[5.2.2]{Gru-Ray} qu'il existe un
\'eclatement $\pi : \ov{X} \to \ov{X}_1$ centr\'e hors de $X$ tel que le transform{\'e} strict
$\ov{\ECC}$ de $\ECC_1$ soit plat. 

$$\xymatrix{ & \ov{X} \ar[d]^\pi\\
X \ar@{^{(}->}[r]_{j_1} \ar@{^{(}->}[ur]^{j} & \ov{X}_1. }$$

Or $\ov{\ECC}$ est coh{\'e}rent car c'est le transform\'e strict d'un faisceau coh\'erent, il est donc localement libre. L'{\'e}clatement ayant lieu
en dehors de $X$, $(\ov{X}, j)$ est encore une compactification de $X$
et $j^*\ov{\ECC} = \ECC.$ Le faisceau $\ov{\ECC}$ est de rang $r$.
\findem

\rqe si $r=1$ {\it i.e,} si $\ECC$ est un faisceau inversible, on peut d\'emontrer cette propri\'et\'e de mani\`ere \'el\'ementaire. Comme pr\'ec\'edemment, on choisit une compactification $(\ov{X}_1, j_1)$ de $X$ et $\ECC_1$ un sous-faisceau coh\'erent de ${j_1}_*(\ECC)$ tel que
$j_1^*(\ECC_1)=\ECC$. D\`es lors on regarde les fibr\'es projectifs :
$$\xymatrix{ \PM(\ECC) \ar@{^{(}->}[r]^j \ar[d]^{\wr} & \PM(\ECC_1) \ar[d] \\ X \ar@{^{(}->}[r]^{j_1} & \ov{X}_1 }$$
On a alors $j^*(\OCC_{\PM(\ECC_1)}(1)) = \OCC_{\PM(\ECC)}(1).$

\begin{theo}
\label{theo9}
Soient $X$, $X'$ deux $k$-vari{\'e}t{\'e}s et un morphisme $f : X' \to X$.
Il existe des compactifications $(\ov{X}, j)$ et $(\ov{X}', j')$ de $X$
et $X'$ respectivement, telles qu'on puisse trouver un morphisme
propre $\ov{f} : \ov{X}' \to \ov{X}$ de sorte que l'on ait le
diagramme commutatif :
$$\xymatrix{X' \ar[d]_f  \ar@{^{(}->}[r]^{j'} & \ov{X}' \ar[d]_{\ov{f}} \\
X \ar@{^{(}->}[r]_{j} & \ov{X}.}$$
De plus, si on se donne un faisceau $\ECC$ localement libre de rang
$r$ sur $X$, on peut choisir ces compactifications de sorte qu'il
existe un faisceau localement libre $\ov{\ECC}$ sur $\ov{X}$ tel que :
$$j^*(\ov{\ECC}) = \ECC.$$
On a alors :
$$(j')^*(\ov{f}^*\ov{\ECC})=f^*\ECC.$$
\end{theo}
\dem On se donne des compactifications $(\ov{X}, j)$ et $(\ov{X}'_1, j')$
de $X$ et $X'$ respectivement. On regarde alors $\ov{X}'_2$
l'adh{\'e}rence sch{\'e}matique de $X'$ dans $\ov{X} \times \ov{X}'_1$. On
note alors $\ov{X}'$ la vari\'et\'e obtenue en \'eclatant le
compl{\'e}mentaire de $X'$ dans $\ov{X}'_2$. L'{\'e}clatement ayant lieu en dehors de $X'$
l'immersion de $X'$ dans $\ov{X}'_2$ se rel{\`e}ve en une immersion de
$X'$ dans $\ov{X}'$. On a donc le diagramme commutatif :
$$\xymatrix{ & \ov{X}' \ar[d]^{\pi} &  \\ X' \ar@{^{(}->}[ur]^{j'} \ar[dd]_f \ar@{^{(}->}[r] & \ov{X}'_2 \ar[dd]_{\ov{f}_2} \ar@{^{(}->}[dr]& \\ & & \ov{X}\times \ov{X}_1' \ar[dl]_{pr}
  \\ X \ar@{^{(}->}[r]_j & \ov{X}. & }$$
La fl\`eche $\ov{f}=\pi \circ \ov{f}_2$ est donc propre comme composition de morphismes propres.
De plus, si on se donne sur $X$ un faisceau localement libre $\ECC$, on
peut gr\^ace au th{\'e}or{\`e}me \ref{theo1} choisir la compactification
$(\ov{X}, j)$ pour qu'il existe sur $\ov{X}$ un faisceau localement
libre $\ov{\ECC}$ v{\'e}rifiant $j^*\ov{\ECC}=\ECC.$
La derni\`ere propri\'et\'e est formelle.
\findem

\begin{theo}
\label{theo2}
 Soient $X$ une $k$-vari{\'e}t{\'e}   
$$0 \to \ECC' \to \ECC \to \ECC'' \to 0,$$
une suite exacte de faisceaux localement libres sur $X$.
Il existe $j : X \hookrightarrow \ov{X}$ une compactification de $X$ et une suite exacte de faisceaux localement libres
$$0 \to \ov{\ECC}' \to \ov{\ECC} \to \ov{\ECC}'' \to 0$$
sur $\ov{X}$, telles que la restriction \`a $X$ de cette derni\`ere soit la suite exacte dont on est parti..
\end{theo}

\dem On choisit une compactification $(\ov{X}_1,j_1)$ de $X$. Il
existe alors un sous-faisceau de ${j_1}_*(\ECC)$ qui soit coh{\'e}rent
et dont la restriction {\`a} $X$ soit $\ECC$. Nommons un tel faisceau
$\ov{\ECC}_1$. On regarde alors la fl{\`e}che compos{\'e}e :
$$\varphi : \ov{\ECC}_1 \to {j_1}_*(\ECC) \to {j_1}_*(\ECC'').$$
On pose :
$${\ov{\ECC}_1}' := \ker \varphi \hbox{ et } {\ov{\ECC}_1}'':=
\ov{\ECC}_1/{\ov{\ECC}_1}'.$$
On obtient donc une suite exacte :
$$0 \to {\ov{\ECC}_1}' \to {\ov{\ECC}_1} \to {\ov{\ECC}_1}'' \to 0$$
o{\`u} ${\ov{\ECC}_1}'$, ${\ov{\ECC}_1}$ et ${\ov{\ECC}_1}''$ sont trois
faisceaux coh{\'e}rents dont les restrictions {\`a} $X$ sont
respectivement $\ECC'$, $\ECC$ et $\ECC''$.


Comme ci-dessus, il existe un {\'e}clatement $p_1 : \ov{X}_2 \to
\ov{X}_1$ en dehors de $X$ tel que $\ov{\ECC}_2 = \TS(\ov{\ECC}_1)$
soit localement libre. En posant :
$${\ov{\ECC}_2}'' := \TS({\ov{\ECC}_1}'') \hbox{ et }{\ov{\ECC}_2}' :=
\ker (\ov{\ECC}_2 \to {\ov{\ECC}_2}'')$$
on a une suite exacte :
$$0 \to {\ov{\ECC}_2}' \to {\ov{\ECC}_2} \to {\ov{\ECC}_2}'' \to 0.$$
De plus ${\ov{\ECC}_2}'$ et ${\ov{\ECC}_2}''$ sont coh{\'e}rents,
${\ov{\ECC}_2}$ est localement libre et la restriction de ces trois
faisceaux {\`a} $X$ est la bonne. On consid{\`e}re alors un {\'e}clatement
$p_2 : \ov{X} \to \ov{X}_2$ en dehors de $X$ (vu comme ouvert de
$\ov{X}_2$) tel que $\ov{\ECC}''= \TS({\ov{\ECC}_2}'')$ soit localement libre. On
proc{\`e}de comme pr{\'e}c{\'e}demment en posant :
$$\ov{\ECC} =  p_2^*({\ov{\ECC}_2})\hbox{ et }  {\ov{\ECC}}' :=
\ker (\ov{\ECC} \to {\ov{\ECC}}'').$$
On obtient alors une suite exacte :
$$0 \to \ov{\ECC}' \to \ov{\ECC} \to \ov{\ECC}'' \to 0$$
o{\`u} les deux faisceaux de droite sont localement libres par
construction (le transform\'e strict d'un faisceau localement libre est localement libre) et celui de gauche comme noyau d'une application surjective entre
faisceaux localement libres. Les restrictions {\`a} $X$ sont toujours
les bonnes.
\findem

\subsection{La premi{\`e}re classe}
\label{standard}
Gr\^ace au th{\'e}or{\`e}me \ref{theo1} on va d{\'e}finir la premi{\`e}re classe
de Chern d'un faisceau inversible $\FCC$ comme l'image par le morphisme de fonctorialit{\'e} de la
premi{\`e}re classe de Chern d'un prolongement de $\FCC$ sur une
compactification. Pour cela, nous avons besoin du th{\'e}or{\`e}me suivant.

\begin{theo}
\label{theo3}
Soient $X$ une vari{\'e}t{\'e} sur $k$ et $\FCC$ un faisceau inversible sur
$X$. Pour $l \in \{1,2\}$ on se donne $(\ov{X}_l, j_l)$ une
compactification de $X$ et $\ov{\FCC}_l$ un faisceau
inversible prolongeant $\FCC$. On a alors :
$${j_{1}}_{rig}^*(c_1(\ov{\FCC}_1)) = {j_{2}}_{rig}^*(c_1(\ov{\FCC}_2)).$$
\end{theo}

\dem

Nous allons d'abord nous ramener {\`a} deux prolongements de notre faisceau
sur la m{\^e}me vari{\'e}t{\'e}. Nous utiliserons pour cela la m{\'e}thode
habituelle du plongement diagonal. On a le diagramme suivant :
$$\xymatrix{
    &  \ov{X}_1  &   &  \\
X \ar@{^{(}->}[ur]^{j_1} \ar[rr]^\alpha \ar@{^{(}->}[dr]_{j_2} &  & \ov{X}_1 \times
\ov{X}_2 \ar[ul]_{p_1} \ar[dl]^{p_2}  & \ov{X_3} \ar[l]_\pi \\
   & \ov{X}_2 &  & }$$

\noindent o{\`u} $\pi$ est l'{\'e}clatement de l'adh{\'e}rence sch{\'e}matique de $X$ dans
$\ov{X}_1 \times \ov{X}_2$ en le compl{\'e}mentaire de $X$. Les
fonctorialit{\'e}s des classes de Chern dans le cas propre nous ram{\`e}nent {\`a}
montrer notre th{\'e}or{\`e}me pour $(p_1 \circ \pi)^*(\ov{\FCC}_1)$ et $(p_2
\circ \pi)^*(\ov{\FCC}_2)$, tous deux faisceaux inversible sur $\ov{X}_3$. En effet si on note encore $\alpha$ l'inclusion de $X$ dans $\ov{X}_3$ on a pour tout $l \in \{1,2\}$ :
\begin{eqnarray}
{j_{l}}_{rig}^*(c_1(\ov{\FCC}_l)) & = & \alpha^* (p_l \circ \pi)^* (c_1(\ov{\FCC}_l)) \\
  & = & \alpha^*c_1((p_l \circ \pi)^*\ov{\FCC}_l).
\end{eqnarray}

De plus, grace \`a la multiplicativit\'e de la premi\`ere classe de Chern, on peut supposer de plus que $\ov{\FCC}_2 = \OCC_{\ov{X}}$ et $\FCC = \OCC_X$.

On suppose donc donn\'ee une compactification $(\ov{X},j)$ et un faisceau inversible $\ov{\FCC}$ tel que 
$$j^*\ov{\FCC} = \OCC_X,$$
et on va montrer que $j_{rig}^*c_1(\ov{\FCC}) = 0.$

Nous allons utiliser un plongement et des calculs {\`a} la Cech pour
faire la preuve. Comme \`a la section \ref{propre}, on se donne une
immersion ferm{\'e}e $\ov{X} \hookrightarrow \YCC$ dans un
$C$-sch{\'e}ma formel lisse au voisinage de $\ov{X}$. On consid{\`e}re
sa fibre g{\'e}n{\'e}rique rigide $\YS$. On a donc le diagramme
commutatif suivant :
$$\xymatrix{X \ar[dr] \ar@{^{(}->}[rr] & & \ov{X}\ar[r]  \ar[dl] & \YCC \ar[d] & \YS
  \ar[l] \ar[d] \\
 & \spec{k} \ar[rr] & & \spf{C} & \spm{K}. \ar[l] }$$

En appliquant le lemme \ref{propre}.\ref{debile} on sait qu'il
existe un recouvrement affine fini $\UG$ de $\YCC$ tel que le recouvrement
induit $\UG_{\ov{X}}$ sur $\ov{X}$ trivialise $\ov{\FCC}$.
On pose alors $\UCC_i =\spf{\ACC_i}$ et $\ov{X}_i = \ov{X} \cap \UCC_i
= \spec{A_i}$ avec $A_i = \ACC_I/I_i$. On pose $X_i = \ov{X}_i \cap X$
et $\UG_{X}$ le recouvrement de $X$ par ces ouverts.
Quitte {\`a} prendre un
recouvrement plus fin, on peut supposer qu'il existe $\ov{f}_i \in A_i$ tel
que $X_i$ soit l'ouvert $D(\ov{f}_i)$ de
$\ov{X}_i$. Pour tout $i$, on choisit alors $f_i \in \ACC_i$ un
rel{\`e}vement de $\ov{f}_i$.

\vspace{.3cm}

\rqe on notera avec un multi-indice $\u{i}=(i_0, \ldots, i_n)$ les
m{\^e}mes choses d{\'e}finies par rapport {\`a} l'ouvert $\UCC_{\u{i}}=
\UCC_{i_0} \cap \cdots \cap \UCC_{i_n}$.

Pour tout multi-indice $\u{i} = (i_0, \ldots, i_n)$, $X_{\u{i}}$ est
l'ouvert $D(\ov{f}_{i_0}.\cdots .\ov{f}_{i_n})$ de
$\ov{X}_{\u{i}}$. On pose donc :
$$\ov{f}_{\u{i}} = \ov{f}_{i_0}.\cdots .\ov{f}_{i_n}.$$
On prend alors comme rel{\`e}vement de $\ov{f}_{\u{i}}$ dans
$\ACC_{\u{i}}$ l'{\'e}l{\'e}ment :
$$f_{\u{i}} := f_{i_0}.\cdots .f_{i_n}$$
o{\`u} $f_{i_k}$ est vu comme {\'e}l{\'e}ment de $\ACC_{\u{i}}$ via la
fl{\`e}che canonique $\ACC_{i_k} \to \ACC_{\u{i}}.$

On choisit alors $u$  un cocycle repr{\'e}sentant la classe
de $\ov{\FCC}$ dans $H^1(\UG_{\ov{X}},
\OCC_{\ov{X}}^*).$

\begin{prop}
Avec les notations pr\'ec\'edentes et en notant $c_1(u)$ l'~image de $c_1(u)$ par l'application $C^2(\UG_K, \Omega^\star_{]\ov{X}[}) \to C^2(\UG_K, j^\dag\Omega^\star_{]\ov{X}[})$, il existe $\zeta \in C^1(\UG_K, j^\dag\Omega^\star_{]\ov{X}[})$ tel que l'on ait l'{\'e}galit{\'e}
$$c_1(u) = \Delta(\zeta).$$
\end{prop}

On sait que ce cocyle repr{\'e}sente un faisceau se restreignant \`a $\OCC_X$ sur $X$. Il existe donc $\theta \in C^0(\UG_X,\OCC_X^*)$ tel que :
$$u = \delta(\theta)$$
o{\`u} $u$ d{\'e}signe par abus de notation l'image dans
$C^1(\UG_X, \OCC_X^*)$ du cocyle $u$.
Pour tout $i$, $\theta_i$ est un \'el\'ement iversible de $(A_i)_{\ov{f}_i}$, il existe donc $k \in \NM$ tel que pour tout $i$ il existe $\alpha_i \in A_i$ et $\beta_i \in A_i$ tels que l'on ait dans ${A_i}_{\ov{f}_i}$ :
$$\frac{\alpha_i}{\ov{f}^k_i}=\theta_i  \hbox{ et } \frac{\beta_i}{\ov{f}^k_i}= \theta^{-1}_i.$$
Quitte \`a augmenter $k$ et \`a modifier $\alpha_i$ et $\beta_i$, on peut supposer que 
$$\alpha_i.\beta_i = \ov{f}^{2k}_i.$$ 

De m\^eme, par d{\'e}finition de $\theta$, on sait que pour tout $i,j \in I$ on a dans $(A_{ij})_{\ov{f}_{ij}}$ :
$$\theta_i u_{ij} = \theta_j.$$
Quitte \`a augmenter de nouveau $k$, on obtient, dans $A_{ij}$, en multipliant par $\ov{f}_{ij}^k$ :
\begin{equation}
\label{lalala2}
\ov{f}_j^k\alpha_i u_{ij} =  \alpha_j \ov{f}_i^k.
\end{equation}
 
On choisit des rel{\`e}vements $\widetilde{\alpha}_i$ (resp. $\widetilde{\beta}_i$) de $\alpha_i$ (resp. $\beta_i$)
dans $\ACC_i$.
On introduit alors pour tout $0 < \lambda < 1$ les voisinages stricts de
$]X_i[$ dans $]\ov{X}_i[$ :
$$V_\lambda^i = \{x\in ]\ov{X}_i[ | |f_i(x)| \geq \lambda \}$$
o\`u on consid\`ere $f_i$ comme une fonction analytique sur $]\ov{X}_i[$.
Sur ces ouverts, la fonction $f_i$ est inversible, on peut donc consid{\'e}rer les fonctions :
$$\widetilde{\theta}_i = \frac{\widetilde{\alpha}_i}{f_i^k}.$$

\begin{lemme} 
Avec les notations pr\'ec\'edentes, il existe un voisinage strict de $]X_i[$ dans $]\ov{X}_i[$ sur lequel $\widetilde{\theta}_i$ est inversible.
\end{lemme}

\dem 
La r\'eduction modulo $I_i$ de $\widetilde{\alpha}_{i}.\widetilde{\beta}_{i} \in \ACC_i$ est $\ov{f}_i^{2k}$. Il existe donc $c_{i} \in I_{i}$ tel que :
$$\widetilde{\alpha}_{i}.\widetilde{\beta}_{i} = f_{i}^{2k} + c_{i}.$$

On utilise alors les voisinages standards. On consid{\`e}re une suite
d'{\'e}l{\'e}ments strictement positifs  $(\lambda_n)$ tendant
inf{\'e}rieurement vers $1$ et une suite $(\eta_n)$ tendant inf{\'e}rieurement vers $1$ telle que pour tout $n$ on ait :
$$\frac{\eta_n}{\lambda_n^{2k}} < 1.$$
On se donne $(g_{i,r})_r$ un syst\`eme de g\'en\'erateurs de l'id\'eal $I_i$. On pose alors :
$$V^i_n = \{x \in ]\ov{X}_i[| |f_i(x)| \geq \lambda_n, \forall r |g_{i,r}(x)|\leq \eta_n\}.$$
On note $V^i_{\u{\lambda},\u{\eta}} = \cup_{n} V^i_n$. C'est un voisinage strict de $]X_i[$ dans $]\ov{X}_i[$ \cite[1.2.4]{Berth5}.

Pour tout $x$ dans cet ouvert on a par d{\'e}finition :
$$\left|\frac{c_{i}}{f_{i}^{2k}}(x)\right| < 1.$$

Donc pour tout $x \in V^i_{\u{\lambda}, \u{\eta}}$ la s{\'e}rie :
$$\sum_{l=0}^\infty (-1)^l\left(\frac{c_{i}}{f_{i}^{2k}}(x)\right)^l$$
converge. Nous la noterons $(1+c_i/f_i^{2k})^{-1}$. La fonction d{\'e}finie sur $V^i_{\u{\lambda},
  \u{\eta}}$ :
$$\widetilde{\theta}_i^{-1} =
\frac{\widetilde{\beta}_i}{f_i^{k}}.\left(1+\frac{c_i}{f_i^{2k}}\right)^{-1}$$
est donc l'inverse de $\widetilde{\theta}_i$.
\findem

On pose alors :
\begin{equation}
\label{zetai}
\zeta_i = \frac{d\widetilde{\theta}_i}{\widetilde{\theta}_i} \in
C^0(\UG_K, j^\dag \Omega^1_{]\ov{X}[}).
\end{equation}

En relevant l'{\'e}quation \ref{lalala2} dans $\ACC_{ij}$, on trouve qu'il existe $b_{ij} \in I_{ij}$ tel que :
$$f_j^k\widetilde{\alpha}_{i}\widetilde{u}_{ij} =\widetilde{\alpha}_jf_i^k + b_{ij}$$
o\`u $\widetilde{u}_{ij}$ est un rel\`evement de $u_{ij}$ servant \`a calculer la classe de Chern. 
On d\'efinit alors les ouverts 
$$V^{ij}_n = \{x \in ]\ov{X}_ij[| |f_{ij}(x)| \geq \lambda_n, \forall r |g_{ij,r}(x)|\leq \eta_n\},$$
o\`u les $g_{ij,r}$ sont des g\'en\'erateurs de $I_{ij}$. On consid\`ere, comme pr\'ec\'edemment $V^{ij}_{\u{\lambda}, \u{\eta}} = \cup_n V^{ij}_n$. C'est un voisinage strict de $]X_{ij}[$ dans $]\ov{X}_{ij}[$.
En voyant tous ces termes comme des fonctions sur $]\ov{X}_{ij}[$ puis en se restreignant aux voisinages stricts $V^{ij}_{\u{\lambda}, \u{\eta}}$, on vient de voir que $\widetilde{\theta}_j$ \'etait inversible.On obtient alors :
\begin{equation}
\label{pourdem}
\widetilde{\theta}_i \widetilde{u}_{ij} = \widetilde{\theta}_j \left( 1 + \frac{b_{ij}}{f_{ij}^k}\widetilde{\theta}_j^{-1}\right).
\end{equation}

De m\^eme, pour tout $x \in V^{ij}_{\u{\lambda}, \u{\eta}},$
$$|\widetilde{\theta}_j^{-1}(x)| \leq \frac{1}{|f_j^k(x)|} \leq \frac{1}{\lambda_n^k}$$
car $|f_i(x)|$ \'etant inf\'erieure \`a $1$, on a 
$$|f_i(x)||f_j(x)| \geq \lambda_n \Rightarrow |f_j(x)| \leq \lambda_n.$$

Donc pour tout $x \in V^{ij}_{\u{\lambda}, \u{\eta}}$ on a :
$$\left| \frac{b_{ij}}{f_{ij}^k}\widetilde{\theta}_j^{-1}(x)\right| < 1.$$
Sur l'ouvert $ V^{ij}_{\u{\lambda}, \u{\eta}}$ on peut donc d{\'e}finir la fonction :
$$\zeta_{ij} := -\log \left(1+\frac{b_{ij}}{f_{ij}^k}\widetilde{\theta}_j^{-1}\right) = \sum_{l=1}^\infty \frac{1}{l}\left(-\frac{b_{ij}}{f_{ij}^k}\widetilde{\theta}_j^{-1}\right)^l.$$

On d{\'e}finit ainsi :
\begin{equation}
\label{zetaij}
\zeta_{ij} \in C^1(\UG_K, j^\dag\OCC_{]\ov{X}[}).
\end{equation} 
En regroupant \ref{zetai} et \ref{zetaij}, on obtient :
$$\zeta \in C^1(\UG_K, j^\dag\Omega^\star_{]\ov{X}[}) = C^1(\UG_K, j^\dag\OCC_{]\ov{X}[}) \oplus C^0(\UG_K, j^\dag \Omega^1_{]\ov{X}[}).$$

Un calcul similaire \`a celui de \ref{bordel}, qui se m\`ene en utilisant l'\'egalit\'e \ref{pourdem}, montre que ce $\zeta$ v\'erifie bien
$$c_1(u) = \Delta(\zeta).$$
\findem

Avec les notations pr{\'e}c{\'e}dentes, $c_1(u)$ est envoy\'e sur z\'ero par l'application :
$$H^2(\UG_K, \Omega^\star_{]\ov{X}[}) \to H^2(\UG_K, j^\dag\Omega^\star_{]\ov{X}[}).$$
On a donc :
$${j}_{rig}^*(c_1(\ov{\FCC})) = 0.$$

\findem

On peut donc maintenant d\'efinir la premi\`ere classe de Chern d'un faisceau inversible.

\begin{Defi}
Avec les notations pr\'ec\'edentes, on pose
$$c_1(\FCC) := j^*c_1(\ov{\FCC})$$
o\`u $(\ov{X},j)$ est une compactification v\'erifiant la condition du lemme \ref{lemme1} et $\ov{\FCC}$ un faisceau inversible sur $\ov{X}$ 
tel que $j^*\ov{\FCC} =\FCC.$
\end{Defi}

\begin{prop}[Fonctorialit{\'e}]
Soient $X$ et $X'$ deux $k$-vari{\'e}t{\'e}s et un morphisme $f : X' \to
X$. Pour tout faisceau inversible $\FCC$ sur $X$ on a :
$$f^*c_1(\FCC) = c_1(f^*\FCC).$$
\end{prop}

\dem Gr\^ace au th{\'e}or{\`e}me \ref{theo9}, on sait qu'il existe des
compactifications $(\ov{X},j)$ et $(\ov{X}',j')$ de $X$ et $X'$
respectivement telles que l'on ait le diagramme commutatif suivant :
$$\xymatrix{X' \ar[d]_f  \ar@{^{(}->}[r]^{j'} & \ov{X}' \ar[d]_{\ov{f}} \\
X \ar@{^{(}->}[r]_{j} & \ov{X}.}$$
De plus, on peut les choisir telles qu'il existe un faisceau inversible
$\ov{\FCC}$  sur $\ov{X}$ v{\'e}rifiant $j^*\ov{\FCC}=\FCC$. On a alors :
$$\begin{array}{ccc} f^*c_1^X(\FCC) & = &
  f^*j^*c_1^{\ov{X}}(\ov{\FCC}) \\
  & = & {j'}^*\ov{f}^*c_1^{\ov{X}}(\ov{\FCC}) \\ & = &
  {j'}^*c_1^{\ov{X}'}(\ov{f}^*\ov{\FCC}) \\ & = & c_1^{X'}(f^*\FCC),
\end{array}$$
la derni{\`e}re {\'e}galit{\'e} venant du fait que l'on a :
$${j'}^*(\ov{f}^*\ov{\FCC}) = f^*\FCC.$$
\findem

\begin{prop}[Multiplicativit\'e]
Soient $X$ une $k$-vari\'et\'e et $\FCC$ et $\FCC'$ deux faisceaux inversibles sur $X$, on a :
$$c_1(\FCC \otimes \FCC') = c_1(\FCC) + c_1(\FCC').$$
\end{prop}

\dem Cela d\'ecoule directement de la multiplicativit\'e de la classe de Chern dans le cas propre \ref{multi}.
\findem

\subsection{Cohomologie de l'espace projectif}

On se pose maintenant la question du calcul de la cohomologie rigide
d'un fibr{\'e} projectif sur une base quelconque. Soient $X$ une
$k$-vari{\'e}t{\'e} et $\ECC$ un faisceau localement libre de rang $r$
sur $X$. On sait gr\^ace au th{\'e}or{\`e}me \ref{theo1} qu'il existe une
compactification $(\ov{X}, j)$ de $X$, ainsi qu'un faisceau localement
libre $\ov{\ECC}$ sur $\ov{X}$ de rang $r$ tel que :
$$j^*\ov{\ECC} = \ECC.$$
On a donc le diagramme commutatif suivant :
$$\xymatrix{\PM(\ECC) \ar@{^{(}->}[r]^{j_P} \ar[d] & \PM(\ov{\ECC}) \ar[d]^p \\
X \ar@{^{(}->}[r]^j & \ov{X} }$$

A partir de l{\`a}, on note $\PM$ pour $\PM(\ECC)$ et $\ov{\PM}$ pour
$\PM(\ov{\ECC})$. On note aussi $\OCC_{\PM}(1)$ et $\OCC_{\ov{\PM}}(1)$
les deux faisceaux canoniques. Si on regarde alors
$c_1(\OCC_{\PM}(1))$ et $c_1(\OCC_{\ov{\PM}}(1))$ comme des
applications dans la cat{\'e}gorie d{\'e}riv{\'e}e $D^+(\ov{X}_{zar})$, on
a le diagramme suivant :
$$\xymatrix{ & &\RM p_* \RM \u{\Gamma}_{rig}(\ov{\PM}/K)[2] \ar[dd] \\
\ZM_{\ov{X}} \ar[urr]^{c_1(\OCC_{\ov{\PM}}(1))}
\ar[drr]_{c_1(\OCC_{\PM}(1))} & & \\
&  & \RM p_*\RM \u{\Gamma}_{rig}(\PM/K)[2]}$$

De m{\^e}me, on a les morphismes de fonctorialit{\'e} suivants :
$$\xymatrix{\RM\u{\Gamma}_{rig}(\ov{X}/K) \ar[r]
  \ar[d] & \RM p_* \RM\u{\Gamma}_{rig}(\ov{\PM}/K)\ar[d] \\
\RM\u{\Gamma}_{rig}(X/K) \ar[r] & \RM p_* \RM\u{\Gamma}_{rig}(\PM/K) }$$

On en d{\'e}duit alors, par cup-produit, comme dans le cas classique
$$\xymatrix{\bigoplus_{i=0}^r c_1(\OCC_{\ov{\PM}}(1))^i : & \RM\u{\Gamma}_{rig}(\ov{X}/K)[-2i] \ar[r] \ar[d] & \RM p_* \RM\u{\Gamma}_{rig}(\ov{\PM}/K) \ar[d] \\
\bigoplus_{i=0}^r c_1(\OCC_{\PM}(1))^i : & 
\RM\u{\Gamma}_{rig}(X/K)[-2i] \ar[r] & \RM p_* \RM\u{\Gamma}_{rig}(\PM/K) }$$

\begin{theo}
\label{theo5}
Avec les notations pr{\'e}c{\'e}dentes, les applications :
$$\bigoplus_{i=0}^r c_1(\OCC_{\ov{\PM}}(1))^i \hspace{1cm} \hbox{et}
\hspace{1cm} \bigoplus_{i=0}^r c_1(\OCC_{\PM}(1))^i,$$
sont des isomorphismes.
\end{theo}

\dem
Le cas des vari\'et\'es propres a d\'ej\`a \'et\'e trait\'e dans \ref{struc}.

Sinon, l'\'enonc\'e \'etant local sur la base, on peut supposer que $\ECC$ est localement libre. Il suffit alors de remarquer que si on se donne une famille cofinale de voisinages stricts de $]X[$ dans $]\ov{X}[$, disons $(V_i)_{i\in I}$, la famille des $\PM^r_{V_i}$ est une famille cofinale de voisinages stricts de $\PM^r_{]X[}$ dans $\PM^r_{]\ov{X}[}$. On utilise alors que le morphisme $p_i : \PM^r_{V_i} \to V_i$ est quasi-compact et quasi-s\'epar\'e et commute donc avec les limites inductives filtrantes \cite[0.1.8]{Berth5}. On est donc ramen\'e au cas trait\'e dans \ref{18}.
\findem

\begin{coro}
\label{lalala}
Avec les notations ci-dessus, on a pour tout $n$ la
d{\'e}composition suivante :
$$H^n_{rig}(\PM/K)=\bigoplus_{i=0}^{r-1}H^{n-2i}_{rig}(X/K).\xi^i.$$

\end{coro}

\subsection{Construction des $c_i$}

On d{\'e}finit les classes de Chern sup{\'e}rieures comme dans le cas
classique. On applique la d{\'e}composition du corollaire \ref{lalala} {\`a}
$\xi^r$ et on obtient :
\begin{equation}
\label{decompog}
\xi^r=\sum_{i=1}^r (-1)^{i+1} c_i(\ECC)\xi^{r-i},
\end{equation}
avec
$$c_i(\ECC) \in H^{2i}_{rig}(X/K).$$

\begin{Defi}
Avec les notations pr{\'e}c{\'e}dentes, pour $0 < i \leq r$ on appelle $i$-{\`e}me classe de
Chern de $\ECC$ la classe $c_i(\ECC)$. On pose de plus 
$$c_0(\ECC) = 1 \hbox{ dans } H^0_{rig}(X/K).$$
\end{Defi} 

\rqe En appliquant ce que l'on vient de voir \`a un fibr\'e inversible $\FCC$, on retrouve notre premi\`ere classe de Chern. En effet la d\'ecomposition \ref{decompog} devient alors :
$$c_{1, rig}(\OCC_{\PM(\FCC)}(1)) = p^*c_1(\FCC).$$
Or la fonctorialit\'e de $c_{1, rig}$ nous donne
$$c_{1, rig}(\OCC_{\PM(\FCC)}(1)) = p^*c_{1, rig}(\FCC).$$
On conclut en utilisant que $p^*$ est injectif.

\begin{prop}[Fonctorialit\'e]
Soient $X$ et $X'$ deux vari{\'e}t{\'e}s propres sur $k$ et $f : X' \to X$
un morphisme. Pour tout $\ECC$ faisceau localement libre de rang $r$ sur $X$ et tout $i$, on a 
$$c_{i}(f^*\ECC) = f^*c_{i}(\ECC).$$
\end{prop}

\dem
On note $\PM' = \PM(f^*(\ECC))$ et $\xi' = c_1(\OCC_{\PM'}(1))$. La fonctorialit\'e de la premi\`ere classe de Chern nous assure que 
$$\xi' = f^* \xi.$$
 En appliquant alors $f^*$ \`a la d\'ecomposition \ref{decompog}, on obtient dans $H^{2r}_{rig}(\PM'/K)$ :
$${\xi'}^r = \sum_{i=1}^r (-1)^i f^*c_i(\ECC){\xi'}^{r-i}.$$
La d\'ecomposition \'etant unique on a :
$$c_i(f^*\ECC) = f^*c_i(\ECC).$$

\subsection{Additivit\'e des classes de Chern}
\begin{theo}
Les classes de Chern d\'efinies pr\'ec\'edemment sont additives.
\end{theo}

\dem Soient $X$ une vari\'et\'e et une suite exacte de faisceaux localement libres 
 
$$0 \to \ECC' \to \ECC \to \ECC'' \to 0.$$
On sait gr\^ace au th\'eor\`eme \ref{theo2} qu'il existe une compactification de $X$, $j : X \hookrightarrow \ov{X}$
et une suite exacte de faisceaux localement libres sur $\ov{X}$ 
$$0 \to \ov{\ECC}' \to \ov{\ECC} \to \ov{\ECC}'' \to 0$$
telles que la restriction \`a $X$ de $\ov{\ECC}$ (resp. $\ov{\ECC'}$, $\ov{\ECC''}$) soit $\ECC$ (resp. $\ECC'$, $\ECC''$).
D\`es lors on a pour tout $i$
$$\begin{array}{ccc} c_i(\ECC) & = & j^*c_i(\ov{\ECC}) \\ 
& = & j^* \left(\sum_{j=0}^i c_i(\ov{\ECC}')c_{i-j}(\ov{\ECC}'')\right) \\
& = & \sum_{j=0}^i c_i({\ECC}')c_{i-j}({\ECC}'').
\end{array}$$

\findem

\subsection{Action du Frobenius}
Soit $X$ une $k$-vari\'et\'e. On note $F$ le frobenius absolu sur $X$. Il induit sur les groupes de cohomologie rigide un frobenius :
$$\Phi : H^i_{rig}(X/K) \to H^i_{rig}(X/K).$$
On a alors 
\begin{prop}
Soit $\ECC$ un faisceau localement libre de rang $r$ sur $X$, pour tout $i \leq r$, on a :
$$\Phi(c_i(\ECC)) = p^ic_i(\ECC).$$
\end{prop}
\dem Par la fonctorialit\'e des classes de Chern, notre proposition se ram\`ene \`a
$$c_i(F^*\ECC) = p^i c_i(\ECC).$$
De plus, de par la d\'efinition des classes de Chern, il suffit de regarder le cas du $c_1$ d'un faisceau inversible. Or dans ce cas, on a $F^*\ECC \iso \ECC^{\otimes p}$. Notre proposition d\'ecoule alors de la multiplicativit\'e des classes de Chern.
\findem

\bibliography{bibthese}

\end{document}